\documentclass[12pt,a4paper,twoside,final,notitlepage, leqno]{article}
\usepackage[english]{babel}
\usepackage[T1]{fontenc}
\usepackage{epsfig, graphicx, amssymb}
\usepackage{amsmath,amsthm,epsfig,amsfonts}
\usepackage{float}
\usepackage{color}
\newcommand{\monitem}{ \smallskip \noindent $\bullet$ \quad  }
\newcommand{\moneq}{\vspace*{-7pt} \begin{equation} \displaystyle }
\newcommand{\moneqstar}{\vspace*{-6pt} \begin{equation*} \displaystyle }
\newcommand{\monendstar}{\vspace*{-6pt} \end{equation*}   }
\newcommand{\monend}{\vspace*{-7pt} \end{equation}   }
\newcommand{\moneqarraystar}{ \begin{eqnarray*} \displaystyle }
\newcommand{\monendarraystar}{ \end{eqnarray*}   }

\newcommand{\dd}{{\rm d}}
\newcommand{\R}{\mathbb{R}}

\definecolor{vertfonce}{rgb}{0.0, 0.5, 0.0}

\def\section*#1{}

\usepackage{fancyhdr}
\renewcommand{\headrulewidth}{0pt}

\begin{document}

\fancypagestyle{plain}{ \fancyfoot{} \renewcommand{\footrulewidth}{0pt}}
\fancypagestyle{plain}{ \fancyhead{} \renewcommand{\headrulewidth}{0pt}}

~

  \vskip 2.1 cm

\centerline {\bf \LARGE General fourth-order Chapman-Enskog expansion}

\bigskip

\centerline {\bf \LARGE of lattice Boltzmann schemes}

 \bigskip  \bigskip \bigskip

\centerline { \large    Fran\c{c}ois Dubois$^{ab}$, Bruce M. Boghosian$^{c}$  and Pierre Lallemand$^{d}$}

\smallskip  \bigskip

\centerline { \it  \small
  $^a$   Laboratoire de Math\'ematiques d'Orsay, Facult\'e des Sciences d'Orsay,}

\centerline { \it  \small   Universit\'e Paris-Saclay, France.}

\centerline { \it  \small
$^b$    Conservatoire National des Arts et M\'etiers, LMSSC laboratory,  Paris, France.}


\centerline { \it  \small
  $^c$ Department of Mathematics, Tufts University, Medford, MA, 02155, USA.} 

\centerline { \it  \small
 $^d$ Beijing Computational Science Research Center, Haidian District, Beijing 100094,  China.}


\bigskip  \bigskip

\centerline {26 August 2023 
  {\footnote {\rm  \small $\,$
    This contribution is published in
{\it Computers and Fluids}, volume 266, article 106036 [11 pages],  November 2023. 
It     has been presented 
at the 31th International Conference on Discrete Simulation of Fluid Dynamics, Suzhou (China) the 22 August 2022.}}}




\bigskip  \bigskip
\noindent {\bf \large Abstract}

\noindent 
In order to derive the equivalent partial differential equations
of a lattice Boltzmann scheme, 
the Chapman Enskog expansion 
is very popular in the lattice Boltzmann community.
A main drawback of this approach is the fact that multiscale expansions
are used without any clear mathematical signification  of the various variables and operators. 
Independently of this framework, 
the Taylor expansion method allows to obtain
formally  the equivalent partial differential equations.
The general equivalency of these two approaches remains an open question.
In this contribution, we prove that both approaches give identical results with 
acoustic scaling for a very general family of lattice Boltzmann schemes
and up to fourth-order accuracy.
Examples with a single scalar conservation illustrate our purpose.

\noindent

\newpage

\bigskip \bigskip    \noindent {\bf \large    1) \quad  Introduction} 

\fancyhead[EC]{\sc{Fran\c{c}ois Dubois, Bruce Boghosian  and Pierre Lallemand}}
\fancyhead[OC]{\sc{General fourth-order Chapman-Enskog expansion}}
\fancyfoot[C]{\oldstylenums{\thepage}}

\smallskip \noindent
The Chapman-Enskog method is a fundamental approach developed for the asymptotic analysis
of the Boltzmann equation. The book of Chapman and Cowling, first published in 1939~\cite{CC39}, contains
the essential of this subject.
When lattice gas automata were first developed in the 1970's~\cite{HPP76} and the 1980's~\cite{FHP86},
the length of the lattice vectors was uniformly equal to unity. The asymptotic analysis for the emergence
of the Navier Stokes equations was conducted by taking the size of the included bodies bigger and bigger.  Employing this limit, a fundamental work was achieved  by H\'enon~\cite{He87} for the determination
of the viscosity of a lattice gas.
With the lattice Boltzmann schemes in the 1990's~\cite{HJ89,HSB89}, the underlying paradigm of the Boltzmann equation
in the approximation of Bhatnagar-Gross-Krook~\cite{BGK54} collision operator became very popular. 
A method of analysis was developed by Chen-Doolen~\cite{CD98} and Qian-Zhou~\cite{QZ00}
based on a Chapman-Enskog expansion. This method was also used by d'Humi\`eres~\cite{DDH92} when he introduced
the multiple relaxation time variant of the lattice Boltzmann schemes.
This approach involves a rather strange formal calculus of partial derivatives with respect to fast and slow time scales that can have noncommutative properties.  Nevertheless, the approach has enjoyed significant success and is recommended in textbooks on the subject, e.g.~\cite{GS13,KKKSSV17,Su01}.

\smallskip \noindent
When one of us began to work in the lattice Boltzmann community, 
lattice Boltzmann schemes were considered as a special finite-difference method on cartesian meshes. 
From this perspective, the classical approach to finding equivalent partial differential equations~\cite{LP73,SY68,WH74}
provided a simple way to make explicit the continuum limit of a given algorithm.
By adapting this method to lattice Boltzmann schemes, we created the Taylor expansion method~\cite{fd07,fd08,fd09,fd22}.
This method has predicted super-convergence of various lattice Boltzmann schemes~\cite{ADGG13,DL09,DL11,DLT10,LD15,LDL23}, 
and elucidated the specific behavior
of the scheme for several sets of boundary conditions~\cite{DLT10,DLT19,DLT20}. 
A natural question is that of precisely how these two approaches, the Chapman-Enskog and Taylor expansion methods, are related. 
This is the subject of this contribution.

\smallskip \noindent
In the second section, we review the Bhatnagar-Gross-Krook framework
and the Chapman-Enskog analysis in the case of a single conserved quantity.
Multi-resolution time lattice Boltzmann schemes are presented in Section 3. 
In the following section, a linear model with a single conserved quantity is presented for two spatial dimensions,
and a preliminary result establishes the equivalence of Chapman Enskog and Taylor  approaches
in this specific case. 
The main result is presented at fourth-order accuracy for very general schemes
in Section 5. The proof for orders three and four is detailed in the two last sections of the paper.

\bigskip \bigskip     \noindent {\bf \large    2) \quad  Bhatnagar-Gross-Krook framework} 

\noindent 
In this section, we follow the standard ``BGK'' framework~\cite{BGK54}
for lattice Boltzmann schemes.
For completeness of our study and to make this work self-contained, we recall known results derived in~\cite{BC98,CD98,fd08,QZ00} relating to Chapman-Enskog and Taylor expansions
for an advection-diffusion model. 

\smallskip \noindent 
At a vertex $ \, x \, $ of a discrete lattice $\, {\cal L} \,$ 
and at discrete time $ \, t $, a  particle distribution with~$ \, q \, $ velocities, 
$ \, f(x,t) = \{f_j(x,t) \, | \,0 \leq j < q\}$, is defined. Its  evolution relative to time  follows a classical algorithm.
First an equilibrium particle distribution $ \, f_j^{\rm eq} (x,t) \, $ is computed from the vector
$ \,  f(x,t) \, $,  
according to a process which is not detailed at this step (see {\it e.g.}~\cite{QDL92}). 
Then the nonlinear relaxation, parametrized by a relaxation time 
$ \, \tau $, is achieved by locally modifying the particle distribution~$ \, f \, $ into a new distribution  $ \, f^* $, defined by the relation
\moneqstar
f_j^*  (x,t) = f_j(x,t) + {{\tau_0}\over{\tau}}\, \big(f_j^{\rm eq} (x,t) -  f_j(x,t) \big) \,,\,\, 0 \leq j < q , 
\monendstar 
with the introduction of a reference time scale $ \, \tau_0 $.
The second step of the algorithm is pure  linear advection of each component of the distribution at its corresponding velocity $\, v_j \, $:
During a small time step  $ \, \Delta t = \varepsilon \, \tau_0 $, the particles stream from the vertex~$ \, x \, $
]to the neighbouring vertex $ \, x + v_j \, \varepsilon \, \tau_0 \, $ of the lattice. 
An iteration of the scheme is  written
%
\moneqstar
f_j(x + v_j \, \varepsilon \, \tau_0 , \, t+\varepsilon \, \tau_0) =  f_j^*  (x , t)  . 
\monendstar 
%
Because $ \, f_j(x + v_j \, \varepsilon \, \tau_0 , \, t+\varepsilon \, \tau_0) =  f_j (x , t) + {{\tau_0}\over{\tau}} \, ( f_j^{\rm eq} -  f_j ) $,
a discrete equation solved by the numerical scheme is easy to make explicit: 
\moneq \label{schema-sm-equilibre} 
f_j(x,\,t)   + {{\tau}\over{\tau_0}} \, \big[ f_j(x + v_j \, \varepsilon \, \tau_0 ,\, t+\varepsilon\, \tau_0)
-  f_j(x,\,t) \big] =  f_j^{\rm eq}  (x , t).
\monend 
With P.\ Coveney and one of us~\cite{BC98},
the linear advection operator 
$ \, D_j \equiv  \partial_t + v_j^\alpha \, \partial_\alpha \, $
in the direction number $ \, j \, $ of the lattice was introduced, 
with an implicit summation on the spatial index $ \, \alpha$. 
Then we can express the  linear advection in terms of the exponential of this operator $\, \exp \, (\varepsilon \,\tau_0 \, D_j) $,
\moneqstar
f_j(x + v_j \, \varepsilon \, \tau_0, \, t+\varepsilon\, \tau_0) = \exp \, (\varepsilon \,\tau_0 \, D_j) \, f_j (x,t) .
\monendstar
With the notation $ \, {\rm I} \, $ for the identity operator and after a second-order expansion of the exponential operator
relative to the small parameter $ \, \varepsilon $,   we
obtain an approximate expression of the functional  equation (\ref{schema-sm-equilibre}): 
\moneq \label{schema-sm-equilibre-approche} 
\Big[ {\rm I} + \varepsilon \, \tau \, \Big( D_j + {{\varepsilon}\over{2}} \, (D_j)^2 + {\rm O}(\varepsilon^2)
\Big) \Big] \, f_j =  f_j^{\rm eq}
\monend
%

\monitem
At this point, the Chapman-Enskog expansion proposed in~\cite{CD98,QZ00} 
introduces  a formal multiple scale expansion for the time derivative,
\moneqstar
\partial_t \equiv \partial_{t_1} + \varepsilon \, \partial_{t_2} + {\rm O}(\varepsilon^2) .
\monendstar
Then the advection operator $ \, D_j \, $ can be expanded in terms of $ \, \varepsilon $: 
\moneqstar
D_j =  D_j^1 + \varepsilon \,  \partial_{t_2}  + {\rm O}(\varepsilon^2) ,\,\, 
D_j^1 = \partial_{t_1}  + v_j^\alpha \, \partial_\alpha . 
\monendstar
We then suppose an {\it a priori} asymptotic expansion  of the particle distribution,
\moneqstar
f \equiv f^0 + \varepsilon \,f^1 + \varepsilon^2 \,f^2 +  {\rm O}(\varepsilon^3),
\monendstar
in the approximate functional equation (\ref{schema-sm-equilibre-approche}) satisfied by the scheme to obtain 
\moneq \label{hierarchie} \left\{ \begin{array}{l}
f_j^0 = f_j^{\rm eq}  \\
f_j^1 + \tau \, D_j^1 f_j^0 = 0   \\
f_j^2 + \tau \, D_j^1 f_j^1 + \tau \, \big[ \partial_{t_2} + {{\tau_0}\over2} \, (D_j^1)^2 \big] f_j^0 = 0 .
\end {array} \right. \monend 
From these relations, we deduce various evolution equations for the distinct time scales $ \, \partial_{t_j} $.

To fix these ideas, we detail the process for  one  conservation law.  In
this case, there is only one scalar conserved variable and we have typically
\moneqstar
\sum_j f_j = \sum_j  f_j^{\rm eq}  \equiv \rho 
\monendstar
with the condition $ \,\, \sum_j  v_j^\alpha \, f_j^{\rm eq} \equiv \rho \, u^\alpha $.
Then  $ \,\, \sum_j f_j^1 = \sum_j f_j^2 = 0 $.
When we insert this condition in the second equation of (\ref{hierarchie}), we obtain
$ \,\, \tau \,  \sum_j (  \partial_{t_1}  + v_j^\alpha \, \partial_\alpha ) f_j^0 = 0  $.
After division by $ \, \tau $, the evolution equation at first order,
\moneq   \label{premier-ordre} 
\partial_{t_1} \rho + u^\alpha \,  \partial_\alpha \rho = 0,
\monend 
is established.

We next insert the condition $ \,\, \sum_j f_j^2 = 0 \,\, $ in the third relation of  (\ref{hierarchie}). 
After dividing by   $ \, \tau $, we obtain
\moneq   \label{second-ordre-preliminaire} 
 \sum_j   D_j^1 f_j^1 +   \sum_j \partial_{t_2} f_j^{\rm eq} +  {{\tau_0}\over2} \, \sum_j (D_j^1)^2 f_j^{\rm eq} = 0 .  
\monend 
We have from the second relation of  (\ref{hierarchie}): 
$ \,\,  \sum_j   D_j^1 f_j^1 = -\tau \,  \sum_j   D_j^1 ( D_j^1 f_j^{\rm eq} ) =  -\tau \,  \sum_j  (D_j^1)^2 f_j^{\rm eq} $,
so the previous relation (\ref{second-ordre-preliminaire})  can be written
\moneq   \label{second-ordre-provisoire} 
\partial_{t_2} \rho + \Big( {{\tau_0}\over2} - \tau \Big) \, \sum_j (D_j^1)^2 \,  f_j^{\rm eq}  = 0,
\monend 
where we also have
%
\begin{align*}
\sum_j (D_j^1)^2 f_j^{\rm eq}
= &
\sum_j
( \partial_{t_1}  + v_j^\alpha \, \partial_\alpha ) \, ( \partial_{t_1}  + v_j^\beta \, \partial_\beta )  f_j^{\rm eq} \\
%
=&
\partial_{t_1}^2 \rho + 2 \, \partial_\alpha \, \partial_{t_1}  (\sum_j   v_j^\alpha \,  f_j^{\rm eq} ) 
+ \sum_j  v_j^\alpha \, v_j^\beta \,  \partial_\alpha \partial_\beta   f_j^{\rm eq} \\
%
=&
\partial_{t_1} ( - u^\alpha \,  \partial_\alpha \rho )
+ 2 \,  \partial_\alpha  \, \partial_{t_1} (  \rho \, u^\alpha ) 
+  \partial_\alpha \partial_\beta (\sum_{j, \alpha,\beta}  v_j^\alpha \, v_j^\beta \,   f_j^{\rm eq} ) \\
%
=&
- u^\alpha \,\partial_\alpha  ( - u^\beta \,\partial_\beta \rho ) 
+ 2 \,  u^\alpha \, \partial_\alpha  ( \partial_{t_1}   \rho ) 
+  \partial_\alpha \partial_\beta (\sum_{j, \alpha,\beta}  v_j^\alpha \, v_j^\beta \,   f_j^{\rm eq} ) \\
%
=&  \partial_\alpha \partial_\beta \big[ \sum_{j, \alpha,\beta}  v_j^\alpha \, v_j^\beta \,   f_j^{\rm eq} -  u^\alpha \,  u^\beta \, \rho \big] 
%
\end{align*}
Inserting this result in the relation (\ref{second-ordre-provisoire}), we obtain the evolution equation
for the second time scale,
\moneq   \label{second-ordre-final} 
 \partial_{t_2} \rho + \Big( {{\tau_0}\over2} - \tau \Big)  \, \partial_\alpha  \partial_\beta \Big[
\sum_j v_j^\alpha \, v_j^\beta \, f_j^{\rm eq} - u^\alpha \, u^\beta \, \rho \Big]  = 0.
\monend 

In this way, with the multiple-time representation  $ \,\,   \partial_t \equiv \partial_{t_1} + \varepsilon \, \partial_{t_2} + {\rm O}(\varepsilon^2) $,
we have finally from the relations (\ref{premier-ordre}) and (\ref{second-ordre-final})
a derivation of the equivalent partial differential equation up to second  order,
\moneq   \label{chapman-enskog-scalaire} 
 \partial_t  \rho + u^\alpha \,  \partial_\alpha \rho - \varepsilon \,  \Big( \tau - {1\over2} \Big) \, \partial_\alpha  \partial_\beta \Big[
\sum_j v_j^\alpha \, v_j^\beta \, f_j^{\rm eq} - u^\alpha \, u^\beta \, \rho \Big]  =   {\rm O}(\varepsilon^2) .
\monend 
%

\monitem
An alternative to the Chapman-Enskog expansion is the Taylor expansion framework
proposed by one of us~\cite{fd08, fd22}. With this paradigm, 
we do not consider multiple time scales
and we do not introduce any  {\it a priori}  asymptotic representation of the particle distribution. 
Using the BGK framework to fix the ideas, we
solve
the approximate functional equation~(\ref{schema-sm-equilibre-approche})
using a formal power series
relative to the small parameter $ \, \varepsilon $,
\begin{align*}
f_j
=&
\big[  {\rm I} - \varepsilon \, \tau \, \big( D_j + {{\varepsilon \, \tau_0}\over{2}} \, (D_j)^2 \big)
 +  (\varepsilon \, \tau)^2 \, (D_j)^2 + {\rm O}(\varepsilon^3) \big] \,  f_j^{\rm eq}  \\
%
=&
\big[  {\rm I}  - \varepsilon \, \tau \, D_j
+ \tau \, \varepsilon^2 \, \big( \tau - {{\tau_0}\over2} \big) \,  (D_j)^2  + {\rm O}(\varepsilon^3) \big] \,  f_j^{\rm eq} \\
%
=&
f_j^{\rm eq}   - \varepsilon \, \tau \, D_j  f_j^{\rm eq}
+ \tau \, \varepsilon^2 \, \big( \tau -  {{\tau_0}\over2}  \big) \,  (D_j)^2  f_j^{\rm eq} + {\rm O}(\varepsilon^3) 
%
%
\end{align*}
For the case of one scalar conserved quantity, we have 
 $ \, \,  \sum_j f_j = \sum_j  f_j^{\rm eq}  \equiv \rho \,\, $ and 
$ \,\, \sum_j  v_j^\alpha \, f_j^{\rm eq} \equiv \rho \, u^\alpha $.  Then after division by $ \, \tau \, \varepsilon $, 
we have the asymptotic relation 
\moneq   \label{taylor provisoire} 
  \sum_j D_j f_j^{\rm eq} - \varepsilon \, \Big( \tau -  {{\tau_0}\over2} \Big) \,  \sum_j  (D_j)^2 f_j^{\rm eq} = {\rm O}(\varepsilon^2) .
\monend
At first order, we have $ \,\,  \sum_j D_j f_j^{\rm eq} = \partial_t \rho + u^\alpha \, \partial_\alpha \rho =  {\rm O}(\varepsilon) $.
To obtain a result at second order, we have  the following calculation:
\begin{align*}
\sum_j (D_j)^2 f_j^{\rm eq}
=&
\sum_j
( \partial_t + v_j^\alpha \, \partial_\alpha ) \, ( \partial_t  + v_j^\beta \, \partial_\beta )  f_j^{\rm eq} \\
%
=&
\partial_t^2 \rho + 2 \, \sum_j v_j^\alpha \, \partial_\alpha \, \partial_t   f_j^{\rm eq}
+ \sum_j v_j^\alpha \,  v_j^\beta \, \partial_\alpha \, \partial_\beta   f_j^{\rm eq} \\
%
=&
\partial_t (- u^\alpha \, \partial_\alpha \rho )  + 2 \,  \partial_\alpha \, \partial_t ( \sum_j v_j^\alpha  f_j^{\rm eq} ) 
+  \partial_\alpha \, \partial_\beta \big[ \sum_j v_j^\alpha \,  v_j^\beta \,  f_j^{\rm eq} \big] + {\rm O}(\varepsilon) \\
%
=&
- u^\alpha \, \partial_\alpha ( - u^\beta \, \partial_\beta \rho )  
 + 2 \,  \partial_\alpha \, \partial_t ( u^\alpha \, \rho)  
+  \partial_\alpha \, \partial_\beta \big[ \sum_j v_j^\alpha \,  v_j^\beta \,  f_j^{\rm eq} \big] + {\rm O}(\varepsilon) \\
%
=&  \partial_\alpha \, \partial_\beta \big[  \sum_j v_j^\alpha \,  v_j^\beta \,  f_j^{\rm eq} - u^\alpha \, u^\beta \, \rho \big] 
+ {\rm O}(\varepsilon) 
\end{align*}
%
We insert this result into the relation (\ref{taylor provisoire}) and we recover exactly the relation  
(\ref{chapman-enskog-scalaire}).

In this section, in the specific case of a scalar equation in the BGK framework, we have established that the Chapman-Enskog methodology and the Taylor expansion method yield exactly the same equivalent
partial differential equation at second-order accuracy.
In the next sections, we generalize this result for an arbitrary number of conservation laws in the framework of multi-resolution times lattice Boltzmann schemes, and we establish agreement up to fourth-order accuracy.

\bigskip \bigskip    \noindent {\bf \large    3) \quad  Multi-resolution times lattice Boltzmann schemes} 

\smallskip \noindent 
The  multi-resolution times paradigm is an extension of the BGK collision operator presented in the previous section.
It has been formalized by d'Humi\`eres~\cite{DDH92}.  As previously, a discrete particle distribution
of $ \, q \, $ velocities  $ \, f(x,t) = \{\, f_j(x,t)\, | \, 0\leq j < q\}$ is defined  with corresponding discrete velocities
$ \,  v_j \in {\cal V} \, $ 
at a vertex $ \, x \, $ 
of
a discrete lattice $ \, \cal{L} \, $ 
and  at discrete time $ \, t $.
One time iteration, 
leading
to the evaluation
of $ \,  f_j(x,t + \Delta t) $, is composed of two steps. 

\smallskip \noindent  
{\it (i)}  Nonlinear relaxation.  
During this step, a local modification of the particle distribution~$ \, f $, denoted by  $ \, f^* $,
is determined.
First an invertible matrix  $ \, M \, $ transforms the particle distribution $ \, f \, $ into moments~$ \, m $.  
We write $\, m = M \, f \, $,
or in terms of components
$ \, \, m_k \equiv \sum_\ell M_{k \ell} \, f_\ell \, $
for $ \, 0 \leq k < q $.
We split this vector into two blocks. The first block~$ \, W \, $ is composed by the conserved quantities
or macroscopic moments, whereas the second block $ \, Y \, $ determines the nonconserved
or microscopic moments,
\moneq \label{mWY} 
m \equiv \begin{pmatrix} W \\ Y \end  {pmatrix} .
\monend 
After relaxation, the conserved moments do not change:  $ \, W^* = W $. 
Secondly, an equilibrium value  $ \, Y^{\rm eq}  \, $ of the nonconserved moments is introduced;
it is a function of the conserved moments,
\moneq \label{Yequilibre}
Y^{\rm eq} = \Phi(W)  .
\monend 
This function $ \, W \longmapsto  \Phi(W) \, $ is required to be regular 
and can be seen as a discrete Gaussian in reference to the Boltzmann equation
for gas dynamics. Nevertheless, it has {\it a priori} no direct algebraic relation with the
Maxwell-Boltzmann distribution
and is constrained 
only
by symmetry considerations.
Knowledge of the equilibrium function $ \, \Phi \, $  is essential for specifying the
multi-resolution times lattice Boltzmann scheme. 
After relaxation the vector of microscopic moments is modified, and a new vector  $ \, Y^* \, $
is created according to 
%
%
%
%
%
\moneq \label{Yetoile} 
  Y^* = Y + S \, (Y^{\rm eq} - Y ) .
\monend 
The  relation (\ref{Yetoile}) introduces a  relaxation matrix $ \, S $.
This is an invertible square matrix, usually chosen as diagonal,
$ \, S = {\rm diag} \, (s_k) $. The relaxation coefficients $ \, s_k \, $
are dimensionless and strictly positive. They are also an
essential specification of the multi-resolution times lattice Boltzmann scheme.
The moments $  \, m^* \, $ after relaxation combine the two families of moments: 
$ \,  m^* = ( W ,\, Y^* )^{\rm t} $.
Then the  particle distribution after relaxation   $ \,   f^* \, $ is simply determined from
the moments  after relaxation: $ \,  f^* = M^{-1} \, m^*  $.
Observe here that if all the relaxation coefficients  $\, s_k \, $ are identical, {\it id est}  
if $ \, s_k \equiv  {{\tau_0}\over{\tau}} \, $ for all indices $ \, k$, then the
multi-resolution times lattice Boltzmann scheme is identical to the BGK variant~\cite{LL00}. 

\smallskip \noindent  
{\it (ii)}   Linear advection. 
This step is identical  to the BGK framework.
Recall that the velocities~$ \, v_j \, $ are chosen in such a way that
after one time step $ \, \Delta t $,
a particle located in $ \, x \in  \cal{L} \, $ arrives at a new vertex of the lattice:
$ \, x +  v_j \,  \Delta t  \in  \cal{L} $. With the previous notation introduced,
we can formulate a compact description of the lattice Boltzmann advection scheme: 
\moneq \label{mrt} 
f_j (x,\, t+\Delta t) =  f_j^*(x - v_j \, \Delta t , \, t)  \,,\,\,  v_j \in {\cal V}
\,,\,\,  x \in  \cal{L} .
\monend
As a final remark, the present paradigm of multiple relaxation schemes
allows one to take into account multiple distributions
of particles. 
Two or more particle distributions can be introduced in practice.
We just observe that the mapping
$ \, j \longmapsto v_j \, $ is not necessarily injective.

\monitem It should be pointed out that important hypotheses have been made for the asymptotic expansions
proposed in this contribution.
First, the discrete function $ \, f(x, \, t) \,$,   for
$ \, x \,  $ a vertex of the lattice and  $\, t \, $  the discrete time,
is assumed to be the  restriction to the lattice of a very regular function 
 denoted in the same way  $ \, f(x, \, t, \Delta t , \, s_k, \, \cdots ) \,$ 
at   a point of the continuous space $ \, x \in  \R^d \, $
and continuous time~$ \, t   $.
The time step    $\, \Delta t \, $ is an infinitesimal quantity,
and this is also the case for the spatial step $ \, \Delta x $.
Additionally, we adopt an acoustic scaling:  The 
numerical velocity  $ \, \lambda \equiv {{\Delta x}\over{\Delta t}} \, $
is supposed fixed as $ \, \Delta x \, $ and~$\, \Delta t \, $ tend to zero.
Last but not least, the relaxation parameters $ \, s_k \, $  are held fixed
when the lattice $ \, {\cal L} \, $ is made finer and finer. 

\bigskip \bigskip
\noindent {\bf \large    4) \quad  A multi-resolution D2Q9 lattice Boltzmann scheme} 

\begin{figure}    [H]  \centering
\centerline {\includegraphics[width=.33\textwidth]{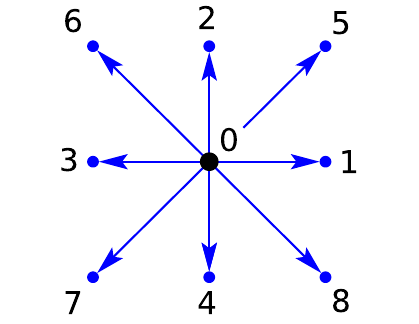}} 

\smallskip  
\caption{D2Q9 lattice Boltzmann scheme }
\label{fig-d2q9} \end{figure}

\smallskip \noindent
In this section, we consider a D2Q9 scheme (see, inter alia,~\cite{CD98,HL97,LL00,QDL92}) 
for a single conservation law in the paradigm of  multi-resolution times lattice Boltzmann schemes. 
This scheme is classic and has been studied in the detail in an article of Luo and one of us~\cite{LL00}. 
The nine velocities begin with $ \, v_0 = 0 \, $ and are presented in the Figure \ref{fig-d2q9}.
The moments $ \, m \, $ are  named as follows in this contribution:
\moneq \label{moments-d2q9}
m^{\rm t} = ( \rho ,\, J_x  ,\,  J_y ,\,  \varepsilon ,\,  XX ,\, XY ,\,  q_x ,\, q_y ,\,  h ) . 
\monend
The density $ \, \rho \, $ is a polynomial of degree zero relative to the velocities,
the momentum~$ \, (J_x  ,\,  J_y ) \, $ is composed by polynomials of degree 1, the energy $ \,  \varepsilon  \, $ and the 
moments $ \, XX \, $ and $ \, XY \, $ are  polynomials of degree 2, the energy flux~$ \, (q_x  ,\,  q_y ) \, $
is associated to polynomials  of degree 3, and the second energy $ \, h \,  $ is of degree~4. 
The explicit construction of the matrix $ \, M \, $ between particles and moments is detailed in the
reference~\cite{LL00}. 
We have 
\moneqstar
M =  \left(     \begin{array} {ccccccccc}
 1  \!&\!  1  \!&\!  1  \!&\!  1  \!&\!  1  \!&\!  1  \!&\!  1  \!&\!  1  \!&\!  1 \\ 
   0 \!&\! \lambda  \!&\!  0  \!&\!  -\lambda  \!&\!  0  \!&\!   \lambda  \!&\!   -\lambda   \!&\!  -\lambda  \!&\!  \lambda  \\ 
 0  \!&\!  0   \!&\! \lambda  \!&\!  0    \!&\!   -\lambda   \!&\!   \lambda   \!&\!   \lambda  \!&\!  -\lambda  \!&\!  -\lambda \\ 
-4 \lambda^2    \!&\!  -\lambda^2   \!&\!   -\lambda^2   \!&\!   -\lambda^2   \!&\!   -\lambda^2  \!&\!   2 \lambda^2   \!&\!  2 \lambda^2   \!&\!  2 \lambda^2   \!&\!  2 \lambda^2 \\ 
  0    \!&\!   \lambda^2    \!&\!  -\lambda^2    \!&\!   \lambda^2   \!&\!   -\lambda^2  \!&\!  0  \!&\!  0  \!&\!  0  \!&\!  0 \\ 
 0  \!&\!  0  \!&\!  0  \!&\!  0  \!&\!  0   \!&\!   \lambda^2   \!&\!   -\lambda^2    \!&\!   \lambda^2   \!&\!   -\lambda^2 \\ 
 0   \!&\! -2 \lambda^3  \!&\!  0   \!&\!  2 \lambda^3  \!&\!  0   \!&\!  \lambda^3   \!&\!   -\lambda^3   \!&\!   -\lambda^3   \!&\!    \lambda^3 \\
 0  \!&\!  0  \!&\! -2 \lambda^3  \!&\!  0  \!&\!  2 \lambda^3    \!&\!   \lambda^3    \!&\!   \lambda^3    \!&\!  -\lambda^3    \!&\!  -\lambda^3 \\ 
4 \lambda^4  \!&\! -2 \lambda^4  \!&\! -2 \lambda^4  \!&\! -2 \lambda^4  \!&\! -2 \lambda^4  \!&\!   \lambda^4   \!&\!  \lambda^4   \!&\!  \lambda^4  \!&\!   \lambda^4 
\end{array} \right) . \monendstar 
%


\newpage 
\monitem 
Advection operator in the basis of moments  

\noindent 
From the velocities $ \, v_j \, $ and the matrix $ \, M $, we introduce the 
momentum-velocity operator matrix (see~\cite{fd22}). It is defined by 
\moneq \label{matrice-Lambda} 
\Lambda  \equiv   M \, {\rm diag} \, \Big(  \sum_\alpha v^\alpha\, \partial_\alpha \Big) \, M^{-1} 
\monend
with  $ \, 1 \leq \alpha  \leq d = $ space dimension. It is simply the set of all advection
operators defined for the lattice and resolved in the basis of the moments. 

\smallskip \noindent
In the case of a single conservation law as studied in the previous section,
there is only one conserved quantity $ \, W \, $ whereas  the vector
$ \, Y \, $ in the relation (\ref{mWY}) is of dimension 8.
We decompose the momentum-velocity operator matrix into four blocks
$\, A $, $ \, B $, $ \, C $, $ \, D \, $ in coherence with the decomposition (\ref{mWY}),
\moneq \label{abcd} 
 \Lambda \equiv  \begin{pmatrix} A &  B \\ C & D \end  {pmatrix} , 
\monend
and similarly for the various powers of  $ \, \Lambda \, $,
\moneqstar 
 \Lambda^2 \equiv  \begin{pmatrix} A_2 &  B_2 \\ C_2 & D_2 \end  {pmatrix} , \,\, 
 \Lambda^3 \equiv  \begin{pmatrix} A_3 &  B_3 \\ C_3 & D_3 \end  {pmatrix} , \,\, 
 \Lambda^4 \equiv  \begin{pmatrix} A_4 &  B_4 \\ C_4 & D_4 \end  {pmatrix} . 
\monendstar
For the D2Q9 scheme and the  advection-diffusion, we have after some lines of algebra
detailed in~\cite{fd22} the following explicit form for the ``ABCD'' decomposition introduced in (\ref{abcd}):

\moneq \label{abcd-d2q9-advection} 
\Lambda_{D2Q9}^{\rm advection}
=
\left( \, \begin{array}{|c|cccccccc|} \hline 
 0 &  \partial_x &  \partial_y &  0  &  0 &  0 &   0 &  0 &  0 
\\  \vspace{-.4 cm} & & & & & & & & \\ \hline  \vspace{-.4 cm} & & & & & & & &  \\
{{2\lambda^2}\over{3}} \, \partial_x &   0 &   0 &  {1\over6} \, \partial_x &
{1\over2} \, \partial_x  &     \partial_y  &    0  &  0  &  0  \\
{{2\lambda^2}\over{3}} \, \partial_y &   0 &   0 &  {1\over6} \, \partial_y & 
 -{1\over2} \, \partial_y  &     \partial_x  &  0  &  0  &  0\\
0 &   \lambda^2 \, \partial_x  &  \lambda^2 \, \partial_y  & 0  &   0  &   0  & 
    \partial_x  &   \partial_y  &    0  \\ 
0  &  {{\lambda^2}\over{3}} \, \partial_x &  -{{\lambda^2}\over{3}} \, \partial_y  & 0 & 0  &  
0  &   - {1\over3} \, \partial_x  &   {1\over3} \, \partial_y   &  0  \\ 
0  &    {{2}\over{3}} \lambda^2 \, \partial_y &    {{2}\over{3}} \lambda^2 \, \partial_x  &  0  & 0  &  
0  &   {1\over3} \, \partial_y  &     {1\over3} \, \partial_x  &   0   \\  
0 &   0 &   0 & {{\lambda^2}\over{3}} \, \partial_x  &  - \lambda^2 \, \partial_x & 
 \lambda^2 \, \partial_y  &    0 &  0  &  {1\over3} \, \partial_x   \\ 
0 &   0 &   0 &  {{\lambda^2}\over{3}} \, \partial_y  &  \lambda^2 \, \partial_y  & 
 \lambda^2 \, \partial_x   &  0  & 0  &    {1\over3} \, \partial_y  \\
 0 &   0 &  0  &  0 &  0 &  0 &  \lambda^2 \, \partial_x  & \lambda^2 \, \partial_y  &  0  \\ \hline 
\end{array} \, \right) .
\monend
%
\vspace{1.5mm}

The structure of the operator matrix $ \, \Lambda \, $ follows ``selection rules'' somewhat similar
to those involved in atoms coupled to the electromagnetic field (see, {\it e.g.}, the book of
Cohen-Tannoudji, Diu and Lalo\"e~\cite{CDL77}). 
The moments
at equilibrium are typically
parametrized  by the two components  $\, u ,\, v \, $  of the imposed velocity and  by a coefficient $ \, \alpha $: 
\moneq \label{Phi-advection-d2q9}
\Phi (\rho) = ( \rho \, u ,\, \rho \, v   ,\, \alpha \, \rho \, \lambda^2   ,\,
 \rho \,  (u^2-v^2) ,\, \rho \, u \, v  ,\, 0  ,\, 0  ,\, 0 \, )^{\rm t} .
\monend
\vspace{1.5mm}
These  moments at equilibrium (\ref{Phi-advection-d2q9}) 
correspond to the following 
nonconserved distribution of particles at equilibrium $ \,  f^{\rm eq}  = M^{-1} \, \Phi  $, with 
\moneqstar \left\{ \begin{array}{l}
f_0^{\rm eq} =  \big( {1\over9} - {{\alpha}\over9} \big) \,  \rho  \\
f_1^{\rm eq} = \big( {1\over9} - {{\alpha}\over36} +  {{u}\over{6\, \lambda}} + {1\over{4\, \lambda^2}} \, (u^2 - v^2) \big) \, \rho \\
f_2^{\rm eq} = \big( {1\over9} - {{\alpha}\over36} +  {{v}\over{6\, \lambda}} - {1\over{4\, \lambda^2}} \, (u^2 - v^2) \big) \, \rho \\
f_3^{\rm eq} = \big( {1\over9} - {{\alpha}\over36} -  {{u}\over{6\, \lambda}} + {1\over{4\, \lambda^2}} \, (u^2 - v^2) \big) \, \rho \\
f_4^{\rm eq} = \big( {1\over9} - {{\alpha}\over36} -  {{v}\over{6\, \lambda}} - {1\over{4\, \lambda^2}} \, (u^2 - v^2) \big) \, \rho \\
f_5^{\rm eq} =   \big( {1\over9} + {{\alpha}\over18}  +  {{u+v}\over{6\, \lambda}}  + {1\over{4\, \lambda^2}} \,u \, v  \big) \, \rho  \\
f_6^{\rm eq} =   \big( {1\over9} + {{\alpha}\over18}  -  {{u-v}\over{6\, \lambda}}  + {1\over{4\, \lambda^2}} \,u \, v  \big) \, \rho  \\
f_7^{\rm eq} =   \big( {1\over9} + {{\alpha}\over18}  -  {{u+v}\over{6\, \lambda}}  - {1\over{4\, \lambda^2}} \,u \, v  \big) \, \rho  \\
f_8^{\rm eq} =   \big( {1\over9} + {{\alpha}\over18}  +  {{u-v}\over{6\, \lambda}}  - {1\over{4\, \lambda^2}} \,u \, v  \big) \, \rho  .
 \end{array} \right. \monendstar
Then we have the relations 
\vspace{2mm}
\moneqstar \left\{ \begin{array}{l}
\sum_j f_j^{\rm eq} = \rho \,,\,\, 
\sum_j v_j^x f_j^{\rm eq} = \rho \, u  \,,\,\,   \sum_j v_j^y f_j^{\rm eq} = \rho \, v  \\ \vspace{-.4 cm} \\ 
\partial_\alpha \partial_\beta \big[ \sum_j v_j^\alpha \, v_j^\beta \, f_j^{\rm eq} - u^\alpha \, u^\beta \, \rho \big]  =
\big( {{\alpha+4}\over6}  - {1\over2} (u^2+v^2)  \big) \, \Delta \rho  
\end{array} \right.
\monendstar

\vspace{1mm}
where $\, \Delta = \partial_x^2 + \partial_y^2 \, $ isthe Laplace operator.
If all the relaxation times
$ \, s_{jx} $, $ \, s_{jy} $, $ \, s_\varepsilon $,  $ \, s_{xx} $,  $ \, s_{xy} $, $ \, s_{qx} $,  $ \, s_{qy} \, $
and $ \, s_h \, $ are equal to the ratio $ \, {{\tau_0}\over{\tau}} $, we have seen in Section~2 that 
the equivalent partial differential equation at second order derived in 
Chapman Enskog derived in (\ref{chapman-enskog-scalaire}) takes the form 
\moneq \label{d2q9-edp-order-2} 
\partial_t  \rho + u \, \partial_x \rho + v \, \partial_y \rho
- \varepsilon \,  \Big( \tau - {{\tau_0}\over2}  \Big) \, 
\Big[ {{\alpha+4}\over{6}} \,\lambda^2  - {1\over2} (u^2+v^2)  \Big] \,  \Delta \rho
 =   {\rm O}(\varepsilon^2) .
\monend
The question now is how to find the equivalent equation
when the relaxation parameters differ.  
Before entering into the resolution of this question, 
we end this section with a general proposition for lattice Boltzmann schemes.

\monitem Exponential expression  of a multi-resolution times lattice Boltzmann scheme

\noindent
We have an exact relation for a discrete time iteration, in the same spirit as for the relation~(\ref{schema-sm-equilibre}). 
It explicitly uses the momentum-velocity operator defined in (\ref{matrice-Lambda}). 

\bigskip \noindent
{\bf Proposition 1: formal expression of one iteration of the scheme}

\noindent
A multi-resolution times lattice Boltzmann scheme (\ref{mrt}) can be written in terms
of  the momentum-velocity operator $\, \Lambda \, $ introduced in (\ref{matrice-Lambda})
through an  exponential operator:
\moneq \label{developpement-exponentiel} 
m  (x, t + \varepsilon \, \tau_0 ) =  {\rm exp} ( - \varepsilon \, \tau_0 \, \Lambda ) \,\,  m^*(x ,\, t) .
\monend 

\smallskip \noindent
The proof of Proposition 1 is given in~\cite{fd22}.  We recall it here to make this contribution self-contained. We have the following calculation: 
\begin{align*}
m_k (x, t + \varepsilon \, \tau_0 )
=&
\sum_{j} M_{k j} \, f_j^* (x-v_j \, \varepsilon \, \tau_0  ,\, t) \\
%
=&
\sum_{j \, \ell} M_{k j} \, (M^{-1})_{_{\scriptstyle \! j \ell}} \, m_\ell^*(x-v_j \, \varepsilon \, \tau_0  ,\, t) \\
%
=&
\sum_{j \, \ell} M_{k j} \, (M^{-1})_{_{\scriptstyle \! j \ell}} \, \sum_{n=0}^{\infty} 
{{1}\over{n !}} \big(-\varepsilon \, \tau_0 \sum_{\alpha}  v_j^\alpha \, \partial_\alpha \big)^n  \, 
m_\ell^*(x ,\, t) \\
%
=&
\sum_{\ell} \sum_{n=0}^{\infty} 
{{1}\over{n !}} \, \sum_{j} M_{k j} \,  \big( \! -\varepsilon \, \tau_0 \, \sum_{\alpha}  v_j^\alpha \, \partial_\alpha \big)^n \, (M^{-1})_{_{\scriptstyle \! j \ell}} 
\, m_\ell^*(x ,\, t) \\
%
=&
\sum_{\ell} \big[  \, \sum_{n=0}^{\infty}  {{1}\over{n !}} \, \big(- \varepsilon \, \tau_0 \, \Lambda \big)^n_{k \ell}  \, \big]  \,\,  m_\ell^*(x ,\, t) \\
%
=&
\sum_{\ell} \,  {\rm exp} ( - \varepsilon \, \tau_0 \, \Lambda )_{k \ell}   \, \,  m_\ell^*(x ,\, t)  \\
%
=&
\big(  {\rm exp} ( - \varepsilon \, \tau_0 \, \Lambda ) \,\, m^*(x ,\, t) \big)_{\! \! k} 
%
\end{align*}
and the relation (\ref{developpement-exponentiel}) is established. 
\hfill $\square $

\bigskip \bigskip    \noindent {\bf \large    5) \quad  Chapman-Enskog framework for  multi-resolution times schemes} 

\smallskip \noindent
In this section, we introduce the Chapman-Enskog expansion in the context of multi-resolution times
lattice Boltzmann schemes. We present at Proposition 2 the multiple times dynamics up to fourth order 
and we deduce general algebraic formulas that control the dynamics at  various scales. We remark  in Proposition 3, that both expansions give identical results.  Then we prove Proposition 2 up to second-order accuracy. 

\smallskip \noindent
In the kinetic theory of gases, 
the small parameter $ \, \varepsilon \, $ in the  Chapman-Enskog expansion is
the ratio of the mean free path, typically 65 nanometers under the usual conditions of temperature and pressure~\cite{Je88}, 
and a characteristic dimension of the problem.
Here, for  multi-resolution times  lattice Boltzmann  schemes, we can set $ \,  \varepsilon = {{\Delta t}\over{\tau_0}}  \, $
as  previously, and the small parameter has a purely numerical interpretation: 
It is the ratio between the time step of the numerical scheme and the reference time scale.
We then expand the particle distribution up to fourth order: 
\moneq \label{dvpt-f-ordre-4}
f = f^{\rm eq} +  \varepsilon \, f^1 + \varepsilon^2 \, f^2 + \varepsilon^3 \, f^3 + {\rm O}(\varepsilon^4 ) .
\monend 

First, we consider an important hypothesis of such Chapman-Enskog expansion: The perturbation terms $ \,   f^\ell  \, $
are  functions only of the equilibrium $ \,  f = f^{\rm eq}  \, $ and its spatial derivatives.
We apply  the d'Humi\`eres matrix $\, M \, $ to the expansion (\ref{dvpt-f-ordre-4}),
\moneqstar
m = M \, f = M \, f^{\rm eq} +  \varepsilon \, M \,  f^1 + \varepsilon^2 \, M \,  f^2 + \varepsilon^3 \, M \,  f^3 +  {\rm O}(\varepsilon^4 ) .
\monendstar
We take the first conserved component of the previous relation. Then $ \,  W = W + 0  \,   $
and the first components of $ \, M \,  f^1 $,  $ \, M \,  f^2 $, {\it etc.} are equal to zero. 

Taking next the  second nonconserved component, we obtain 
\moneqstar
Y = Y^{\rm eq} +   \varepsilon \, (M \,  f^1)_Y  +   \varepsilon^2 \, (M \,  f^2)_Y +   \varepsilon^3 \, (M \,  f^3)_Y  + {\rm O}(\varepsilon^4) \,
\monendstar
and the perturbation terms $ \, \varepsilon^\ell \, (M \,  f^\ell)_Y \, $  depend only
 on  the conserved moments $ \, W \, $ and their spatial derivatives. 
We introduce the specific notations $\,\Psi_j \, $  for the previous expansion: 
\moneq \label{dvpt-Y-ordre-3}
 Y = \Phi(W) + S^{-1} \, \big( \varepsilon \, \tau_0 \, \Psi_1(W) + \varepsilon^2 \, \tau_0^2 \,\Psi_2(W)
 + \varepsilon^3 \, \tau_0^3 \,\Psi_3(W) \big)
+ {\rm O}(\varepsilon^4)  , 
\monend
with $ \, Y^{\rm eq} = \Phi(W) $, see  (\ref{Yequilibre}). 
We suppose also as in~\cite{CD98, QZ00} a 
multi-scale  approach  for the  time dynamics:
\moneqstar
\partial_t =  \partial_{t_1}   + \varepsilon \,  \partial_{t_2}
+ \varepsilon^2 \, \partial_{t_3}  + \varepsilon^3 \, \partial_{t_4} + {\rm O}(\varepsilon^4)  . 
\monendstar
%
 
\bigskip  \noindent
{\bf Proposition 2: multiple time dynamics with the Chapman-Enskog expansion}

\noindent
With the hypotheses presented previously, 
the conserved quantities $ \, W \, $ follow a multiple time dynamics :
\moneq \label{dvpt-dtW-ordre-4}
\partial_{t_1} W \!+ \Gamma_1(W) = 0 \,,\,\, 
\partial_{t_2} W \!+ \tau_0 \, \Gamma_2 (W) = 0  ,\,\, 
\partial_{t_3} W \!+ \tau_0^2 \, \Gamma_3(W) = 0  ,\,\,  
\partial_{t_4} W \!+ \tau_0^3 \, \Gamma_4(W) = 0   .  
\monend
The differential operators
 $ \, \Gamma_1 (W) $, $\, \Psi_1 (W) $, $\, \Gamma_2 (W) $, $\, \Psi_2 (W) $,
$\, \Gamma_3 (W) $, $\, \Psi_3 (W) \, $ and $ \, \Gamma_4 (W) \, $
introduced in the relations (\ref{dvpt-Y-ordre-3}) and (\ref{dvpt-dtW-ordre-4}) 
are  determined recursively as  functions
of the data $ \, v_j $, $ \, M $, $ \, \Phi(W) \, $ and $ \,S $.
The  operator $\,\Gamma_1(W) \, $ establishes the first-order dynamics,
%
%
\moneq \label{ordre-1}
\Gamma_1 = A \, W + B \, \Phi(W).
\monend
After introducing the H\'enon matrix 
\moneq \label{henon}
\Sigma \equiv  S^{-1} - {1\over2} \, {\rm I}
\monend
that generalizes the expansion first presented by H\'enon in~\cite{He87},
the differential operators  $\,\Psi_1(W) \, $ and $\,\Gamma_2(W) \, $
have to be specified for the second-order evolution: 
\moneq \label{ordre-2} \left\{ \begin{array}{l}
\Psi_1 = \dd \Phi(W) .  \Gamma_1  - (C \, W + D \, \Phi(W)) \\ \\
\Gamma_2 = B \, \Sigma \, \Psi_1 .
\end{array} \right. \monend
%
%
%
At  third order, we have 
\moneq \label{ordre-3} \left\{ \begin{array}{l}
\Psi_2  (W) = \Sigma \, \dd \Psi_1 .  \Gamma_1 +  \dd \Phi(W) .  \Gamma_2
-   D \, \Sigma \, \Psi_1  \\ \\
\Gamma_3  (W) = B \, \Sigma \, \Psi_2   - {1\over6} \,B \, \dd \Psi_1  .  \Gamma_1  
+ {1\over12} \, B_2 \, \Psi_1,
\end{array} \right. \monend
%
%
and at   fourth   order
\moneq \label{ordre-4} \left\{ \begin{array}{l}
 \Psi_3  (W) =    \Sigma \, \dd \Psi_1  .  \Gamma_2  +  \dd \Phi .  \Gamma_3 -  D \, \Sigma \, \Psi_2 + \Sigma \, \dd \Psi_2 .  \Gamma_1 
 +{1\over6} \, D \, \dd \Psi_1 .  \Gamma_1 \\   \qquad \qquad 
  - {1\over12} \, D_2 \, \Psi_1 (W)  - {1\over12} \, \dd \, (\dd \Psi_1 .  \Gamma_1 ) .  \Gamma_1  \\ \\
\Gamma_4   (W) = B \, \Sigma \, \Psi_3  + {1\over4} \, B_2 \, \Psi_2   +  {1\over6} \, B \, D_2 \, \Sigma \, \Psi_1
-   {1\over6} \, A \, B \, \Psi_2 -  {1\over6} \, B \, \dd \, (\dd \Phi .  \Gamma_1)  . \Gamma_2 \\   \qquad \qquad 
  -  {1\over6} \, B  \, \dd \, (\dd \Phi .  \Gamma_2) . \Gamma_1   - {1\over6} \, B \, \Sigma \,\dd \, (\dd \Psi_1 .  \Gamma_1 ) .  \Gamma_1. 
\end{array} \right. \monend


\monitem 
The proof of this proposition constitutes the remainder of this contribution.
We first observe, however, that with the Taylor expansion method
an asymptotic partial differential system  is emerging~\cite{fd22}: 
\moneq \label{edp-taylor-ordre-4}
\displaystyle \partial_t W + \Gamma_1 + \Delta t \, \Gamma_2  +  \Delta t^2 \, \Gamma_3
      +  \Delta t^3 \, \Gamma_4  = {\rm O}(\Delta t^4) .
\monend
The coefficients $\, \Gamma_j \, $ in the expansion (\ref{edp-taylor-ordre-4})
are  vectors  obtained after $ \, j \, $ spatial derivations  of the conserved moments $ \, W \, $ 
and the equilibrium vector $ \, \Phi(W) $.  For the non-conserved moments, we have 
\moneq \label{taylor-Y-ordre-3}
Y = \Phi(W) +   S^{-1} \, \big( \Delta t \, \Psi_1
  + \Delta t^2 \, \Psi_2  + \Delta t^3 \, \Psi_3    \big) +  {\rm O}(\Delta t^4)  . 
\monend
The  differential operators $\, \Psi_j \ $ are analogous to $\, \Gamma_j \, $ but not with the same dimension.
The explicit forms of the operators 
 $ \, \Gamma_1 (W) $, $\, \Psi_1 (W) $, $\, \Gamma_2 (W) $, $\, \Psi_2 (W) $,
$\, \Gamma_3 (W) $, $\, \Psi_3 (W) \, $ and $ \, \Gamma_4 (W) \, $
relative to the fourth-order Taylor expansion have been derived in our contribution~\cite{fd22}. 
We have the following result. 

\bigskip \noindent
{\bf Proposition 3: The Taylor and Chapman-Enskog expansions give identical results}

\noindent
With the hypotheses presented at the end of Section 3, the 
precise
algebraic expression of the operators 
$ \, \Gamma_1 (W) $, $\, \Psi_1 (W) $, $\, \Gamma_2 (W) $, $\, \Psi_2 (W) $,
$\, \Gamma_3 (W) $, $\, \Psi_3 (W) \, $ and $ \, \Gamma_4 (W) \, $
are identical to the relations
(\ref{ordre-1}), (\ref{ordre-2}), (\ref{ordre-3}) and (\ref{ordre-4}).

\smallskip \noindent
The proof of this proposition is
obtained by comparing
the results of Proposition 2
with
the main result of our previous contribution. 
We just have to observe that the expressions~(\ref{edp-taylor-ordre-4})
and (\ref{taylor-Y-ordre-3}) use the same notations as in the reference~\cite{fd22}, and that the relations (\ref{ordre-1}) to (\ref{ordre-4}) are exactly the same as
those proposed in this reference. 
\hfill $\square$ 

\monitem Example: advection-diffusion with the D2Q9 scheme

\noindent
Before entering into the different steps of the proof of Proposition 2, 
we illustrate the previous expansion with the  scalar conservation law
studied in the previous section.
The moments are still given by the relation (\ref{moments-d2q9})
and the nonconserved moments at equilibrium
by the formulas~(\ref{Phi-advection-d2q9}). 
%
The operator matrix $ \, \Lambda_{D2Q9}^{\rm advection} \, $ for  advection-diffusion  has been made explicit
in~(\ref{abcd-d2q9-advection}). 
The block decomposition (\ref{abcd}) can be determined. We have $\, A = 0 \, $ and
\moneq \label{B-d2q9}
B \, (j_x ,\, j_y ,\, \varepsilon ,\, xx ,\, yy,\,  q_x ,\, q_y ,\, h)^{\rm t} = \partial_x j_x + \partial_y j_y .
\monend
Then at first order, we have
\moneqstar
\Gamma_1  = A \, W + B \, \Phi(W)  = u \, \partial_x \rho + v \, \partial_y \rho .
\monendstar
Also, we have from (\ref{ordre-1}) that $\, \Psi_1  = \dd \Phi(W) .  \Gamma_1  - (C \, W + D \, \Phi(W))  \, $
and
\moneqstar
(\Psi_1)_{jx} =  {{u^2 + v^2}\over2}  \, \partial_x \rho   - \Big( {2\over3} + {{\alpha}\over6} \Big) \, \lambda^2 \, \partial_x \rho  \,,\,\, 
(\Psi_1)_{jy} =  {{u^2 + v^2}\over2}  \, \partial_y \rho   - \Big( {2\over3} + {{\alpha}\over6} \Big) \, \lambda^2 \, \partial_y \rho .
\monendstar
The H\'enon matrix
$ \,\,   \Sigma \equiv  S^{-1} - {1\over2} \, {\rm I}  \,\, $ is a diagonal matrix and
we impose isotropy conditions: $ \, \sigma_{jx} = \sigma_{jy} = \sigma_j \, $
and $ \, \sigma_{qx} =  \sigma_{qy} = \sigma_q $.  We then have
\moneqstar
\Sigma = {\rm diag} \, \big( \sigma_j ,\, \sigma_j ,\, \sigma_e ,\, \sigma_x ,\, \sigma_x ,\,
\sigma_q ,\, \sigma_q ,\, \sigma_h \big)
\monendstar
with $ \,\,  \sigma_j = {1\over{s_j}} - {1\over2}  $.
Finally, at second order, we have $\,\,  \Gamma_2 = B \, \Sigma \, \Psi_1 $,
and due to the structure~(\ref{B-d2q9}) of the $ \, B \, $ differential operator, only the two
first components $\, (\Psi_1)_{jx} \, $ and $ \, (\Psi_1)_{jy} \, $ 
of the vector $ \, \Psi_1 \, $ are used. Then we have
$ \, B \, \Sigma \, \Psi_1 = \sigma_j \, \big( \partial_x  (\Psi_1)_{jx}  + \partial_y  (\Psi_1)_{jy}  \big) \, $ and
\moneqstar
\Gamma_2  = \sigma_j \, \Big[  {{u^2 + v^2}\over2}  - \Big( {2\over3} + {{\alpha}\over6} \Big) \, \lambda^2 \Big] \, \Delta \rho  .
\monendstar
Finally, due to (\ref{dvpt-dtW-ordre-4}), the equivalent partial differential equation of the D2Q9 lattice Boltzmann scheme is written 
%
%
\moneqstar
\partial_t  \rho + u \, \partial_x \rho + v \, \partial_y \rho
- \varepsilon \, \tau_0 \,  \Big( {1\over{s_j}} - {1\over2} \Big) \, 
\Big[ {{\alpha+4}\over6} \,\lambda^2  - {1\over2} (u^2+v^2)  \Big] \,  \Delta \rho =   {\rm O}(\varepsilon^2) .
\monendstar
This equation is very  similar to the equation (\ref{d2q9-edp-order-2})
established previously for the BGK variant.  Now, we know that with multi-resolution  times lattice Boltzmann schemes,
the coefficient of dissipation is exactly related to the relaxation coefficient $ \, s_j \, $
for  the momentum $\, J $. 

\monitem
Chapman-Enskog expansion: Study at order zero 

\noindent
We establish here that in the expansion (\ref{dvpt-Y-ordre-3}), the first term is the equilibrium
function $ \, Y^{\rm eq} = \Phi(W) $. 
We start from the formal expansion (\ref{developpement-exponentiel}):
$ \,  m (t+\varepsilon \, \tau_0) = \exp ( - \varepsilon \, \tau_0 \,  \Lambda ) \,  m^* $.
At order zero, we can write
\moneqstar
 m +  {\rm O}(\varepsilon) =  m^*  +  {\rm O}(\varepsilon) .
\monendstar
For the first component, we have: $ \,\,  W  +  {\rm O}(\varepsilon) = W^*  +  {\rm O}(\varepsilon) $.
This relation provides no new information because $ \, W^* = W $.
For the second component we obtain
$ \,  Y  +  {\rm O}(\varepsilon) = Y^*  +  {\rm O}(\varepsilon) $. Due to the 
relaxation (\ref{Yetoile}), we have   $ \, Y^* = Y + S \, ( \Phi(W) - Y) $.
The matrix  $ \, S \, $ is supposed fixed and invertible.
Then  
\moneq \label{Y-Ystar-ordre-0}
Y =  \Phi(W)  + {\rm O}(\varepsilon) \,,\,\,
 \, Y^* =  \Phi(W)  + {\rm O}(\varepsilon) . 
\monend
%

\monitem
Chapman-Enskog expansion : Study at order one 

\smallskip \noindent
We consider the expansion (\ref{developpement-exponentiel})
at order one with
$ \, \partial_t =  \partial_{t_1}  + {\rm O}(\varepsilon)  $. Then 
\moneq \label{m-de-t-ordre-1}
m + \varepsilon  \, \tau_0 \,  \partial_{t_1} m   +  {\rm O}(\varepsilon^2)
  =  m^*  - \varepsilon  \, \tau_0 \,  \Lambda \,  m^* +  {\rm O}(\varepsilon^2)  
\monend
with \quad 
$ m^{\rm t}  =  ( W \,,\,\, Y ) \, $ and  $\, \Lambda \, $ decomposed into four blocks according to  (\ref{abcd}).
For the first component of the relation (\ref{m-de-t-ordre-1}), we have 
\moneqstar
W   + \varepsilon  \, \tau_0 \,  \partial_{t_1} W   +  {\rm O}(\varepsilon^2)  = W^*
- \varepsilon  \, \tau_0 \, ( A \, W + B \, Y^* ) +  {\rm O}(\varepsilon^2)  
\monendstar
with  $ \, W^* = W \, $  and   $ \, Y^* = \Phi(W) +  {\rm O}(\varepsilon) $. 
Then $ \,  \partial_{t_1} W = - (  A \, W + B \,  \Phi(W) ) \, $
and the first relation of (\ref{dvpt-dtW-ordre-4})
is established, with  $ \, \Gamma_1 (W) \, $ given by the relation (\ref{ordre-1}). 
\hfill $\square$ 

\monitem Chapman-Enskog expansion: End of the study at order one 

\noindent 
We look now at the second  component of the relation  (\ref{m-de-t-ordre-1}): 
\moneqstar
Y + \varepsilon  \, \tau_0 \,  \partial_{t_1} Y   +  {\rm O}(\varepsilon^2)
  =  Y^*  - \varepsilon  \, \tau_0 \,  (C \, W + D \, Y^* ) +  {\rm O}(\varepsilon^2)  .
\monendstar
Then
$ \,\, Y -  Y^* = -  \varepsilon  \, \tau_0 \, \big( \partial_{t_1} Y +  (C \, W + D \, Y^* ) \big)  +  {\rm O}(\varepsilon^2)  $.
If we take also into consideration the relation (\ref{Yetoile}), we have the exact relation 
$ \,\, S \, ( Y -  \Phi(W)) =  Y - Y^* $.
Then, after taking into consideration the  expansions $ \, Y = \Phi(W) + {\rm O}(\varepsilon) \, $
and  $ \, Y^* = \Phi(W) + {\rm O}(\varepsilon) $, 
we have the following calculation:
\begin{align*}
S \, \big( Y - \Phi(W) \big)
=&
Y - Y^*  \\
%
=&
- \varepsilon \, \tau_0 \,  \partial_{t_1} \big( \Phi(W) +  {\rm O}(\varepsilon) \big)
- \varepsilon  \, \tau_0  \,  (C \, W + D \, \big(  \Phi(W) +  {\rm O}(\varepsilon) \big) +  {\rm O}(\varepsilon^2)  \\
%
=&
 \varepsilon   \, \tau_0 \, \big[ - \dd \Phi(W)  .  \partial_{t_1} W -  (C \, W + D \,  \Phi(W) ) \big]
  + {\rm O}(\varepsilon^2)    \\
%
=&
 \varepsilon  \, \tau_0  \,  \big[ \dd \Phi(W)  . \Gamma_1  -  (C \, W + D \,  \Phi(W) ) \big] 
  + {\rm O}(\varepsilon^2)     \\
%
=&
\varepsilon  \, \tau_0  \, \Psi_1 (W) + {\rm O}(\varepsilon^2),     
%
\end{align*}
with    $\,  \Psi_1 = \dd \Phi(W) .  \Gamma_1  - (C \, W + D \, \Phi(W)), $ and the first relation
of (\ref{ordre-2}) is established.
\hfill $\square$

\monitem
Taking into account the H\'enon matrix inside the expansion 

\smallskip  \noindent
From the relations  (\ref{Yetoile}) and (\ref{dvpt-Y-ordre-3}), we have the two expansions at first order 
\moneqstar
Y = \Phi(W) + \varepsilon  \, \tau_0 \,  S^{-1} \, \Psi_1(W)  + {\rm O}(\varepsilon^2 ) \,,\,\,
Y^* = \Phi(W) + \varepsilon  \, \tau_0 \,  ( S^{-1} - {\rm I} ) \, \Psi_1(W)  + {\rm O}(\varepsilon^2 ) .
\monendstar
With the matrix $\, \Sigma \, $ introduced in (\ref{henon}), 
we have also the expansions at first order 
\moneq \label{Y-Ystar-ordre-1} \left\{ \begin{array}{l}
Y = \Phi(W) + \varepsilon \, \tau_0 \,  \big( \Sigma +  {1\over2} \,{\rm I}  \big) \,\, \Psi_1(W)  + {\rm O}(\varepsilon^2 ) \\ 
Y^* = \Phi(W) + \varepsilon \, \tau_0 \,  \big( \Sigma -  {1\over2} \,{\rm I}  \big) \,\, \Psi_1(W) + {\rm O}(\varepsilon^2 ) .
\end{array} \right. \monend 
%

\monitem Chapman-Enskog expansion : Study at order two

\smallskip \noindent
We consider again the expansion (\ref{developpement-exponentiel}),
but now at order two, and we obtain 
\moneqstar
m + \varepsilon  \, \tau_0 \, \partial_t m   +   {1\over2} \, \varepsilon^2  \, \tau_0^2 \, \partial_t^2 m + {\rm O}(\varepsilon^3) =
m^*  - \varepsilon  \, \tau_0 \,  \Lambda \,  m^*  +   {1\over2} \, \varepsilon^2  \, \tau_0^2 \, \Lambda^2 \,  m^*  +  {\rm O}(\varepsilon^3) .
\monendstar
We introduce the multiple scales for time evolution:
$ \, \partial_t =  \partial_{t_1}  + \varepsilon \, \partial_{t_2}   + {\rm O}(\varepsilon^2)  $. Then

\smallskip \noindent 
$ m + \varepsilon  \, \tau_0 \,  ( \partial_{t_1} +   \varepsilon \, \partial_{t_2} ) \,  m   +   {1\over2} \, \varepsilon^2  \, \tau_0^2 \,
 ( \partial_{t_1} +   {\rm O}(\varepsilon) )^2  m  + {\rm O}(\varepsilon^3)  
 = m^*  - \varepsilon  \, \tau_0 \,  \Lambda \,  m^*  +   {1\over2} \, \varepsilon^2  \, \tau_0^2 \, \Lambda^2 \,  m^*  +  {\rm O}(\varepsilon^3) $

\smallskip \noindent 
and we have
\moneq \label{m-de-t-ordre-2}
m + \varepsilon  \, \tau_0 \, \partial_{t_1} m + \varepsilon^2  \, \tau_0 \, \Big( \partial_{t_2}   m  +   {{\tau_0}\over2} \, \partial_{t_1}^2 m \Big)  
 = m^*  - \varepsilon  \, \tau_0 \,  \Lambda \,  m^*  +   {1\over2} \, \varepsilon^2  \, \tau_0^2 \, \Lambda^2 \,  m^*  +  {\rm O}(\varepsilon^3) .
\monend
The square of the operator  $ \,   \Lambda \, $ satisfies 
 $ \,   \Lambda^2  = \begin{pmatrix} A & B \\ C & D \end{pmatrix} \begin{pmatrix} A & B \\ C & D \end{pmatrix}
\equiv \begin{pmatrix} A_2 & B_2 \\ C_2 & D_2 \end{pmatrix}  \, $
and 
\moneq \label{A2-B2-C2-D2}
  A_2  = A^2 + B \, C \,,\,\,   B_2  = A \, B + B \, D \,,\,\, C_2 = C \, A + D \, C  \,,\,\, D_2 = C \, B + D^2  .
\monend
and similar operators for higher powers of the matrix $ \, \Lambda $. 
Then the first component of the relation (\ref{m-de-t-ordre-2}) can be written
\moneqstar  \left\{ \begin{array}{l}
W + \varepsilon  \, \tau_0 \, \partial_{t_1} W + \varepsilon^2  \, \tau_0 \, \big( \partial_{t_2}   W
+  {{\tau_0}\over2} \, \partial_{t_1}^2 W \big)  \\ 
\qquad \qquad  =  W  - \varepsilon  \, \tau_0 \,  (A \, W + B \,  Y^*  ) +   {1\over2} \, \varepsilon^2  \, \tau_0^2 \,
( A_2 \, W  + B_2 \, Y^*  )  +  {\rm O}(\varepsilon^3) .
\end{array} \right. \monendstar 
The terms at order zero of the previous relation are eliminated.
At order one, we have to take into account the relation
$ \, \,  Y^*  = \Phi(W) + \varepsilon  \, \tau_0 \,  (\Sigma \, \Psi_1  -  {1\over2} \Psi_1 ) + {\rm O}(\varepsilon^2 )  $.
Then we recover the relation 
$ \, \partial_{t_1} W + A \, W  + B \, \Phi(W) = 0 \, $ established previously.
At second order a new relation is emerging: 
\moneq \label{dt2-W-provisoire}
 \partial_{t_2}  W +   {{\tau_0}\over2} \, \partial_{t_1}^2 W = - \, \tau_0 \, B \,  \Big( \Sigma \, \Psi_1  -  {1\over2} \, \Psi_1 \Big)  +
 {{\tau_0}\over2} \, \Big( A_2 \, W  + B_2 \, \Phi \Big) .
\monend
From the relation (\ref{ordre-1}), we have
\begin{align*}
\partial_{t_1}^2 W
=&
\partial_{t_1} \, ( -\Gamma_1(W) ) \\
%
=&
 -  \partial_{t_1} \, ( A \, W + B \, \Phi(W) ) \\
%
=&
A \,  \Gamma_1 + B \, \dd  \Phi(W) .  \Gamma_1 \\
%
=&
A \, ( A \, W + B \, \Phi) + B \, \dd  \Phi(W) . \Gamma_1  \\
%
=&
A_2  \, W - B\, C \, W + A \, B \, \Phi  + B \, \dd  \Phi(W) . \Gamma_1,  
%
\end{align*}
due to (\ref{A2-B2-C2-D2}). Then the relation (\ref{dt2-W-provisoire}) can be written 
\begin{align*}
\partial_{t_2} W& + \tau_0 \, B \,  \Sigma \, \Psi_1\\
=&
- {{\tau_0}\over2}  \, \big( A_2 W - B\, C \, W + A \, B \, \Phi
+ B \, \dd  \Phi(W) . \Gamma_1 \big)
+   {{\tau_0}\over2} \, B \, \Psi_1 +   {{\tau_0}\over2} \, (   A_2  W  +    B_2 \, \Phi ) \\
%
%
=& \, 
 {{\tau_0}\over2}  \, \big (    B\, C \, W - A \, B \, \Phi - B \, \dd  \Phi(W) . \Gamma_1 
+ B \, ( \dd  \Phi(W) . \Gamma_1 - C \, W - D \, \Phi) + (A \, B + B \, D ) \, \Phi \\
%
=& \, 0 .
%
\end{align*}
This last relation expresses exactly that  
$ \,   \partial_{t_2} W + \tau_0 \, \Gamma_2 (W) = 0  \,  $ with  $ \,  \Gamma_2 (W) =  B \,  \Sigma \, \Psi_1 (W), $
and the second relation of (\ref{ordre-2}) is established. 
\hfill $\square$

\bigskip \bigskip    \noindent {\bf \large    6) \quad  Chapman-Enskog expansion at order three} 

\smallskip \noindent
The relations established in the previous section are very useful.
For example, we have used them to study the ability to recover formally  the compressible Navier-Stokes 
equations at second order with only one particle distribution~\cite{DL22}.
In order to study finer properties of the lattice Boltzmann scheme, however, a higher precision is necessary.
We have done this for specific problems in previous contributions~\cite {DL11, DL22, LD15, LDL23}.
Here, we establish general formulas (\ref{ordre-3}) for  future works.  

\monitem Chapman-Enskog expansion: End of the study at order two

\noindent
We first look to the second-order expansion  (\ref{m-de-t-ordre-2}). 
The  second component can be written
\moneqstar \left\{ \begin{array}{l}
Y + \varepsilon \, \tau_0 \,  \partial_{t_1} Y +   \varepsilon^2 \,  \tau_0 \, \partial_{t_2} Y
+   {1\over2} \, \varepsilon^2 \, \tau_0^2  \, \partial_{t_1}^2  Y
\\ \qquad \qquad 
= Y^*  - \varepsilon \, \tau_0 \,  (C \, W + D \, Y^*)
+   {1\over2} \, \varepsilon^2 \, \tau_0^2 \, ( C_2 \, W + D_2 \, Y^* ) +  {\rm O}(\varepsilon^3) .
 \end{array} \right. \monendstar
Then we have

\noindent $
S \, ( Y - \Phi(W) )   =  Y - Y^*  $

\smallskip \noindent $ \qquad \qquad   \qquad \,\,\, =
- \varepsilon \, \tau_0 \, \partial_{t_1}  Y  -  \varepsilon^2 \, \tau_0 \, \big(\partial_{t_2}  Y 
+   {{\tau_0}\over2} \, \partial_{t_1}^2   Y  \big)
- \varepsilon \, \tau_0 \,  (C \, W + D \, Y^*) $

\smallskip \noindent  \qquad \qquad   \qquad   \qquad 
$ +   {1\over2} \, \varepsilon^2 \, \tau_0^2 \, ( C_2 \, W + D_2 \, Y^* ) +  {\rm O}(\varepsilon^3) . $

\smallskip \noindent
We insert the representations (\ref{Y-Ystar-ordre-1}) 
into the right-hand side of the previous expansion to obtain

\noindent $
S \, ( Y - \Phi(W) )   =   - \varepsilon \, \tau_0 \,  \partial_{t_1} \Phi(W)
-  \varepsilon^2 \, \tau_0 \, \big[ \tau_0 \,  \partial_{t_1}  ( \Sigma \, \Psi_1 +   {{\tau_0}\over2} \, \Psi_1 ) +  \partial_{t_2} \Phi(W)
 +  {{\tau_0}\over2} \, \partial_{t_1}^2   \Phi(W) \big] $

\smallskip  \qquad  \qquad  \qquad  \qquad $   
 - \varepsilon \, \tau_0 \,  \big[ C \, W + D \, \big( \Phi(W) + \varepsilon \, \tau_0 \, \big( \Sigma \, \Psi_1 -  {1\over2} \, \Psi_1 \big) \big]
 +   {1\over2} \, \varepsilon^2\, \tau_0^2  \, ( C_2 \, W + D_2 \, \Phi  ) +  {\rm O}(\varepsilon^3) . $ 

\smallskip \noindent
We have by definition 
$  \,\, S \, ( Y - \Phi(W) )  =  \varepsilon \, \tau_0 \, \Psi_1 +  \varepsilon^2 \, \tau_0^2 \,  \Psi_2 +  {\rm O}(\varepsilon^3) $.  
The  first-order term relative to $ \varepsilon \, $  is proportional
to  $ \,  \dd \Phi . \Gamma_1 - (C \, W + D \, \Phi(W) ), \, $ and we recover  $ \, \Psi_1 \, $ 
due to the first relation of~(\ref{ordre-2}). We can make explicit the 
second-order term from the previous calculation: 

\smallskip \noindent
$  \Psi_2 = -  \Sigma \,   \partial_{t_1}  \Psi_1  - {1\over2} \, \partial_{t_1}  \Psi_1 
-  {{1}\over{\tau_0}} \, \partial_{t_2} \Phi(W) -  {1\over2} \, \partial_{t_1}^2   \Phi(W)
- D \,   \big( \Sigma \, \Psi_1 -  {1\over2} \, \Psi_1 \big)
+  {1\over2} \,  C_2 \, W  +  {1\over2}   D_2   \,  \Phi(W)   $ 

\smallskip \noindent
with



\smallskip \noindent
$  \partial_{t_1}  \Psi_1  = \partial_{t_1} \big( \dd \Phi . \Gamma_1 - C \, W - D \, \Phi(W) \big)  $

\smallskip  \qquad  $ \,\,\,\,   =
\partial_{t_1}  (\dd \Phi . \Gamma_1 ) -  C \, \partial_{t_1} W - D \, \dd \Phi .  \partial_{t_1} W $

\smallskip  \qquad  $ \,\,\,\,   =
\partial_{t_1}  (\dd \Phi . \Gamma_1 ) + C \, \Gamma_1 +  D \, \dd \Phi . \Gamma_1  \,, $

\smallskip \noindent
$ {1\over{\tau_0}} \,  \partial_{t_2} \Phi(W)  =  {1\over{\tau_0}} \, \dd  \Phi(W) .  \partial_{t_2} W = -  \dd \Phi(W) . \Gamma_2  \,, $ 

\smallskip \noindent
$ \partial_{t_1}^2   \Phi(W)  = \partial_{t_1} \big(  \partial_{t_1}   \Phi(W) \big)
= \partial_{t_1} ( \dd \Phi .  \partial_{t_1} W ) = - \partial_{t_1} (  \dd \Phi . \Gamma_1 )  \,, $ 

\smallskip \noindent
$ C_2  = C \, A + B \, D  \,$ and $ \,  D_2  = C \, B + D^2 $. 
We deduce

\smallskip  \noindent   
$  \Psi_2 = \Sigma \, \dd \Psi_1 . \Gamma_1
- {1\over2} \, \big( \partial_{t_1}  (\dd \Phi . \Gamma_1 ) + C \, \Gamma_1 +  D \, \dd \Phi . \Gamma_1 \big) 
 + \dd  \Phi(W) . \Gamma_2
+ {1\over2} \, \partial_{t_1} (  \dd \Phi . \Gamma_1 ) 
$ 

\smallskip  \qquad  \quad $ - D \,  \Sigma \, \Psi_1 +  {1\over2} \, D \, \Psi_1  +  {1\over2} \, \ C \, ( A \, W + B \, \Phi)
+  {1\over2} \, D \, ( C \, W + B \, \Phi  ) $

\smallskip  \noindent   \quad $ \,\, =
\Sigma \, \dd \Psi_1 . \Gamma_1   - {1\over2} \,   D \, \dd \Phi . \Gamma_1  +  \dd \Phi(W) . \Gamma_2
- D \,  \Sigma \, \Psi_1 +  {1\over2} \,    D \, \Psi_1
+  {1\over2} \, D \, ( \dd \Phi . \Gamma_1   -   \Psi_1  ) $ 

\smallskip  \noindent   \quad $ \,\, =
\Sigma \, \dd \Psi_1 . \Gamma_1  +  \dd \Phi(W) . \Gamma_2 - D \,  \Sigma \, \Psi_1   $

\smallskip  \noindent
and the first relation of (\ref{ordre-3}) is proven. 
\hfill $\square   $ 

\monitem Chapman-Enskog expansion: Study at order three 

\smallskip \noindent
We refer to the expansion (\ref{developpement-exponentiel}) 
at order three, and we obtain
\moneqstar \left\{ \begin{array}{l}
m + \varepsilon  \, \tau_0 \, \partial_t  m   +   {1\over2} \, \varepsilon^2  \, \tau_0^2 \,   \partial_t^2  m
 +   {1\over6} \, \varepsilon^3  \, \tau_0^3 \,  \partial_t^3  m \\ \qquad \qquad 
 = m^*  - \varepsilon  \, \tau_0 \,  \Lambda \,  m^*  +   {1\over2} \, \varepsilon^2  \, \tau_0^2 \, \Lambda^2 \,  m^*
-   {1\over6} \, \varepsilon^3  \, \tau_0^3 \, \Lambda^3 \,  m^*  +  {\rm O}(\varepsilon^4)  
\end{array} \right. \monendstar 
and 
\moneqstar \left\{ \begin{array}{l}
 m + \varepsilon  \, \tau_0 \, ( \partial_{t_1} +   \varepsilon \, \partial_{t_2}  +   \varepsilon^2 \, \partial_{t_3})  \,  m
+   {1\over2}  \, \tau_0^2 \, \varepsilon^2 \,   ( \partial_{t_1} +  \varepsilon \, \partial_{t_2} + {\rm O}(\varepsilon^2) )^2 \,   m
+   {1\over6} \, \varepsilon^3  \, \tau_0^3 \,  ( \partial_{t_1} + {\rm O}(\varepsilon) )^3  \,  m  
\\ \qquad \qquad = m^*  - \varepsilon  \, \tau_0 \,  \Lambda \,  m^*
+ {1\over2} \, \varepsilon^2  \, \tau_0^2 \, \Lambda^2 \,  m^* - {1\over6} \, \varepsilon^3  \, \tau_0^3 \, \Lambda^3 \,  m^*  +  {\rm O}(\varepsilon^4) .
\end{array} \right. \monendstar 
We expand the various powers of
$\, \partial_t = \partial_{t_1} +   \varepsilon \, \partial_{t_2}  +   \varepsilon^2 \, \partial_{t_3} + {\rm O}(\varepsilon^3) $, 
paying attention to the non commutation of these operators. For example,
$ \, \partial_{t_1} \, \partial_{t_2} \not = \partial_{t_2}  \,  \partial_{t_1} $: 
\moneq \label{m-de-t-ordre-3} \left\{ \begin{array}{l}
m + \varepsilon  \, \tau_0 \,  ( \partial_{t_1} +   \varepsilon \, \partial_{t_2}  +   \varepsilon^2 \, \partial_{t_3}) \, m
 +   {1\over2} \, \varepsilon^2  \, \tau_0^2 \, ( \partial_{t_1}^2 +  \varepsilon \, \partial_{t_1} \, \partial_{t_2} +
 \varepsilon \, \partial_{t_2} \, \partial_{t_1} )  \, m \\ \qquad \qquad
  +   {1\over6} \, \varepsilon^3  \, \tau_0^3 \, \partial_{t_1}^3  \, m  
 = m^*  - \varepsilon \, \tau_0 \,  \Lambda \,  m^*
+ {1\over2} \, \varepsilon^2 \, \tau_0^2 \, \Lambda^2 \,  m^* - {1\over6} \, \varepsilon^3 \, \tau_0^3 \, \Lambda^3 \,  m^*  +  {\rm O}(\varepsilon^4) .
\end{array} \right. \monend
We consider the first component of the relation (\ref{m-de-t-ordre-3}), relative to the conserved variables: 

\smallskip  \noindent 
$ W + \varepsilon \, \tau_0 \,   ( \partial_{t_1} +   \varepsilon \, \partial_{t_2}  +   \varepsilon^2 \, \partial_{t_3})  W
 +   {1\over2} \, \varepsilon^2 \, \tau_0 \, ( \partial_{t_1}^2 +  \varepsilon \, \partial_{t_1} \, \partial_{t_2} +
 \varepsilon \, \partial_{t_2} \, \partial_{t_1} ) W
 +   {1\over6} \, \varepsilon^3 \, \tau_0^3 \, \partial_{t_1}^3 W  $ 

\smallskip  \noindent \qquad 
$ = W  - \varepsilon \,  \tau_0 \, ( A \, W + B \,  Y^*  )
+   {1\over2} \, \varepsilon^2 \, \tau_0^2 \,  ( A_2 \, W + B_2 \,  Y^*  )
 - {1\over6} \, \varepsilon^3 \, \tau_0^3 \,   ( A_3 \, W + B_3 \,  Y^*  )  +  {\rm O}(\varepsilon^4) ,  $
\smallskip  \noindent 
with  $ \,\,     Y^*  = \Phi(W) + \varepsilon \, \tau_0 \,  (\Sigma \, \Psi_1  -  {1\over2} \Psi_1 )
 + \varepsilon^2 \,\tau_0^2 \,  (\Sigma \, \Psi_2  -  {1\over2} \Psi_2 ) + {\rm O}(\varepsilon^3 )  $. 
Then we obtain 
\moneqstar \left\{ \! \begin{array}{l}
W + \varepsilon \, \tau_0 \,  ( \partial_{t_1} +   \varepsilon \, \tau_0^2 \, \partial_{t_2}  +   \varepsilon^2  \, \partial_{t_3})\,  W
 +   {1\over2} \, \varepsilon^2 \,  \tau_0^2 \, ( \partial_{t_1}^2 +  \varepsilon  \, \partial_{t_1} \, \partial_{t_2} +
\varepsilon  \, \partial_{t_2} \, \partial_{t_1} ) \, W 
+   {1\over6} \,  \varepsilon^3  \, \tau_0^3 \, \partial_{t_1}^3  \, W \\ \quad = W  - \varepsilon \, \tau_0 \,  A \, W
- \varepsilon \, \tau_0 \, B \, \big[ \Phi(W) + \varepsilon \,\tau_0 \,   (\Sigma \, \Psi_1  -  {1\over2} \Psi_1 )
 +   \varepsilon^2  \, \tau_0^2 \,  (\Sigma \, \Psi_2  -  {1\over2} \Psi_2 ) \big] \\ \qquad 
 +   {1\over2} \, \varepsilon^2 \, \tau_0^2 \,   A_2 \, W   
 + {1\over2} \, \varepsilon^2 \,  \tau_0^2 \,  B_2 \,  \big[ \Phi(W) + \varepsilon  \, \tau_0 \,  (\Sigma \, \Psi_1  -  {1\over2} \Psi_1 ) \big]
 \\ \qquad
 - {1\over6} \, \varepsilon^3 \,  \tau_0^3 \,  ( A_3 \, W + B_3 \,  \Phi )   +  {\rm O}(\varepsilon^4) .
\end{array} \right. \monendstar 
We identify the third-order terms of the previous relation:
\moneqstar \left\{ \begin{array}{l}
{1\over{\tau_0^2}} \,  \partial_{t_3}  W +  {{1}\over{2\, \tau_0}} \, ( \partial_{t_1} \partial_{t_2} W  +   \partial_{t_2}  \partial_{t_1}  W  )
+ {{1}\over6} \, \partial_{t_1}^3 W  
\\ \quad
= - B \,  (\Sigma \, \Psi_2  -  {1\over2} \Psi_2 )  +
{1\over2} \,  B_2 \, (\Sigma \, \Psi_1  -  {1\over2} \Psi_1 )  - {1\over6} \, ( A_3 \, W + B_3 \,  \Phi )  
\end{array} \right. \monendstar 
with

\smallskip  \noindent
$ {1\over{\tau_0}} \,\partial_{t_1} \partial_{t_2} W  =  \partial_{t_1} ( - B \, \Sigma \, \Psi_1 ) =
-  B \, \Sigma \, \dd \Psi_1 . \partial_{t_1} W =  B \, \Sigma \, \dd \Psi_1 . \Gamma_1   , $ 

\smallskip  \noindent
$   {1\over{\tau_0}} \, \partial_{t_2}  \partial_{t_1}  W =  {1\over{\tau_0}} \, \partial_{t_2} ( - A \, W - B \, \Phi )
= A \, \Gamma_2 -  {1\over{\tau_0}} \, B \, \dd \Phi . \partial_{t_2} W
=   A \, B \, \Sigma \, \Gamma_1 + B \, \dd \Phi .  \Gamma_2   $. 

\smallskip \noindent
Then 

 \noindent
$   {1\over{\tau_0^2}} \,  \partial_{t_3}  W +  {1\over2} \, B \, \Sigma \, \dd \Psi_1 . \Gamma_1 
 +  {1\over2} \, (  A \, B \, \Sigma \, \Gamma_1 + B \, \dd \Phi .  \Gamma_2 ) 
+ {1\over6} \,   \partial_{t_1}^3 W  $

\smallskip  \qquad   $  \,\,\, = 
- B \,  (\Sigma \, \Psi_2  -  {1\over2} \Psi_2 )  +
{1\over2} \,  B_2 \, (\Sigma \, \Psi_1  -  {1\over2} \Psi_1 )  - {1\over6} \, ( A_3 \, W + B_3 \,  \Phi ) $. 

\smallskip \noindent
We observe that 

\smallskip \noindent
$ \partial_{t_1}^3 W  = 
\partial_{t_1} ( A \, \Gamma_1 + B \, \dd \Phi . \Gamma_1 )
= \partial_{t_1} \big( A \, ( A \, W + B \, \Phi ) + B \, \dd \Phi . \Gamma_1 \big)  $

\smallskip  \qquad   $  \,\,\, = - A \, ( A \, \Gamma_1 - B \, \dd \Phi . \Gamma_1 ) -  B \, \dd \, (\dd \Phi . \Gamma_1).\Gamma_1 $

\smallskip  \qquad   $  \,\,\, =  - A^2 \, \Gamma_1 - A \, B \, \dd \Phi . \Gamma_1  -  B \, \dd \, (\dd \Phi . \Gamma_1).\Gamma_1 $, 

\smallskip  \noindent
$ A_3  = A_2 \, A + B_2 \, C \,,\,\, B_3 = A_2 \, B + B_2 \, D $,

\smallskip  \noindent
$ A_3 \, W + B_3 \,  \Phi  = A_2 \, ( A \, W + B \, \Phi) + B_2 \, ( C \, W + D \, \Phi) $ 

\smallskip  \noindent  \qquad  \qquad  \qquad $  =
(A^2 + B\, C) \, \Gamma_1 + (A \, B + B\, D) \, ( \dd \Phi . \Gamma_1 - \Psi_1 ) $ 

\smallskip  \noindent  \qquad  \qquad  \qquad $  =
A \, ( A \, \Gamma_1 + B \, \dd \Phi . \Gamma_1 ) + B \, (C \, \Gamma_1 + D \, \dd \Phi . \Gamma_1 )
- B_2 \, \Psi_1 $ 

\smallskip  \noindent  \qquad  \qquad  \qquad $  =
A \, ( A \, \Gamma_1 + B \, \dd \Phi . \Gamma_1 ) + B \, (   \dd (\dd \Phi . \Gamma_1).\Gamma_1
- \dd \Psi_1 . \Gamma_1 ) - B_2 \, \Psi_1 $ 

\smallskip  \noindent  \qquad  \qquad  \qquad $  =
- \partial_{t_1}^3 W  -  B \,  \dd \Psi_1 . \Gamma_1   - B_2 \, \Psi_1 , $ 

\smallskip  \noindent
$ \Psi_2  = \Sigma \, \dd \Psi_1 . \Gamma_1  +  \dd \Phi .  \Gamma_2 - D \,  \Sigma \, \Psi_1 . $

\noindent
In consequence, we have

\smallskip  \noindent $
 {1\over{\tau_0^2}} \,  \partial_{t_3}  W = -  {1\over2} \,  ( B \, \Sigma \, \dd \Psi_1 . \Gamma_1
+  A \, B \, \Sigma \, \Gamma_1 + B \, \dd \Phi .  \Gamma_2 )
- B \,  \Sigma \, \Psi_2  +   {1\over2}  \,  B \,  \Psi_2  $ 

\smallskip  \noindent \qquad \qquad \quad $
+ {1\over2} \,  B_2 \, (\Sigma \, \Psi_1  -  {1\over2} \Psi_1 )
+ {1\over6} \, (  B \,  \dd \Psi_1 . \Gamma_1  + B_2 \, \Psi_1 )  $

\smallskip  \noindent \qquad  \quad   $ \,\,\,\, =
-  {1\over2} \,  B \, \Sigma \, \dd \Psi_1 . \Gamma_1
-  {1\over2} \,   A \, B \, \Sigma \, \Gamma_1 -  {1\over2} \, B \, \dd \Phi .  \Gamma_2  - B \,  \Sigma \, \Psi_2 $ 

\smallskip  \noindent \qquad \qquad \quad $
+    {1\over2}  \,  B \, (   \Sigma \, \dd \Psi_1 . \Gamma_1 +    \dd \Phi . \Gamma_2 
-  D \,  \Sigma \, \Psi_1  ) 
+ {1\over2} \,   B_2 \, \Sigma \, \Psi_1  - \big(  {1\over4} -  {1\over6} \big) \, B_2 \,  \Psi_1   + {1\over6} \,   B \,  \dd \Psi_1 . \Gamma_1   $

\smallskip  \noindent \qquad  \quad   $ \,\,\,\, =
-  B \,  \Sigma \, \Psi_2 - {1\over12}  \,  B_2 \,  \Psi_1 + {1\over6} \,   B \,  \dd \Psi_1 . \Gamma_1 
$

because $ \, B_2 = A\, B + B \, D $.
Then the third relation of (\ref{dvpt-dtW-ordre-4})  is established, and
$ \, \Gamma_3 \, $ is given by the second relation of (\ref{ordre-3}).
\hfill $\square   $

\bigskip \bigskip    \noindent {\bf \large    7) \quad  Chapman-Enskog expansion at order four} 

\smallskip \noindent
We  establish the first relation of (\ref{ordre-4}) and make explicit the expression
for $ \, \Psi_3 $. Then we extract the value of $ \, \Gamma_4 \, $ from (\ref{developpement-exponentiel}) 
and establish the second  relation of (\ref{ordre-4}).

\monitem Chapman-Enskog expansion: End of the study at order three

\smallskip \noindent
We consider the second component of the relation (\ref{m-de-t-ordre-3}):



 \smallskip  \noindent 
 $ Y + \varepsilon \,\tau_0 \,   ( \partial_{t_1} \, +   \varepsilon  \, \partial_{t_2}  +   \varepsilon^2  \, \partial_{t_3})  Y
 +   {1\over2} \, \varepsilon^2 \,\tau_0^2 \, ( \partial_{t_1}^2 +  \varepsilon \, \partial_{t_1} \, \partial_{t_2} +
 \varepsilon  \, \partial_{t_2} \, \partial_{t_1} ) Y
 +   {1\over6} \,   \varepsilon^3 \,\tau_0^3 \,\,  \partial_{t_1}^3 Y   $ 

 \smallskip  \noindent \qquad 
$ =  Y^* 
 - \varepsilon \, \tau_0 \,  ( C \, W + D \,  Y^*  )  +   {1\over2} \, \varepsilon^2 \,\tau_0^2 \,  ( C_2 \, W + D_2 \,  Y^*  )
 - {1\over6} \,   \varepsilon^3 \,\tau_0^3 \,   ( C_3 \, W + D_3 \,  Y^*  )  +  {\rm O}(\varepsilon^4)  . $

 \smallskip  \noindent
We insert in this relation the representation at order two of  the  nonconserved variables  $ \, Y \, $ and~$ \, Y^* \, $
\moneqstar \left\{ \begin{array}{l}
Y= \Phi(W) + \varepsilon \, \tau_0 \, (\Sigma \, \Psi_1  +  {1\over2} \Psi_1 )
+ \varepsilon^2 \,\tau_0^2 \,  (\Sigma \, \Psi_2  +  {1\over2} \Psi_2 ) + {\rm O}(\varepsilon^3) \\
Y^*  = \Phi(W) + \varepsilon \,\tau_0 \,  (\Sigma \, \Psi_1  -  {1\over2} \Psi_1 )
+ \varepsilon^2 \,\tau_0^2 \,  (\Sigma \, \Psi_2  -  {1\over2} \Psi_2 ) + {\rm O}(\varepsilon^3)   .
\end{array} \right. \monendstar 
Then

\smallskip \noindent
$  Y  +   \varepsilon \,\tau_0 \,  \partial_{t_1} \big[   \Phi(W) + \varepsilon \,\tau_0 \,  (\Sigma \, \Psi_1  +  {1\over2} \Psi_1 )
  + \varepsilon^2 \,\tau_0^2 \,  (\Sigma \, \Psi_2  +  {1\over2} \Psi_2)  \big] $

\smallskip \noindent  \quad 
$ +  \, \varepsilon^2 \, \tau_0 \,  \partial_{t_2} \big[ \Phi(W) + \varepsilon \,\tau_0 \,  (\Sigma \, \Psi_1  +  {1\over2} \Psi_1 ) \big] 
 +  \varepsilon^3 \,\tau_0 \,  \partial_{t_3} \Phi(W) $

\smallskip \noindent  \quad 
$ +   {1\over2} \, \varepsilon^2 \,\tau_0^2 \, \partial_{t_1}^2 \big[ \Phi(W) + \varepsilon \,\tau_0 \,  (\Sigma \, \Psi_1  +  {1\over2} \Psi_1 ) \big] 
+    {1\over2} \, \varepsilon^3 \,\tau_0^2 \,\partial_{t_1} \, \partial_{t_2}  \Phi(W)
+   {1\over2} \, \varepsilon^3 \,\tau_0^2 \, \partial_{t_2} \, \partial_{t_1}  \Phi(W) $

\smallskip \noindent  \quad $
+   {1\over6} \,   \varepsilon^3  \,\tau_0^3 \, \partial_{t_1}^3  \Phi(W)
=  Y^*  - \varepsilon \,\tau_0 \,    C \, W $

\smallskip \noindent  \quad $
- \varepsilon \,\tau_0 \,  D \,  \big[ \Phi(W)   + \varepsilon \,\tau_0 \,   (\Sigma \, \Psi_1  -  {1\over2} \Psi_1 )
+   \varepsilon^2 \,\tau_0^2 \,   (\Sigma \, \Psi_2  -  {1\over2} \Psi_2 ) \big] 
+   {1\over2} \, \varepsilon^2 \,\tau_0^2 \,  C_2 \, W $

\smallskip \noindent  \quad $
+   {1\over2} \,   \varepsilon^2  \,\tau_0^2 \, D_2 \,  \big[ \Phi(W)  +    \varepsilon \,\tau_0 \,  (\Sigma \, \Psi_1  -  {1\over2} \Psi_1 ) \big]
- {1\over6} \,  \varepsilon^3  \,\tau_0^3 \,  ( C_3 \, W  + D_3 \, \Phi )  +  {\rm O}(\varepsilon^4)  $ 

\smallskip \noindent 
We identify the terms relative to   $ \, \varepsilon^3  \, $
in the relation
\moneqstar
Y - Y^*  = S \, (Y - \Phi) =
 \varepsilon \,\tau_0 \, \Psi_1 + \varepsilon^2 \,\tau_0^2 \, \Psi_2 +  \varepsilon^3  \,\tau_0^3 \, \Psi_3  +  {\rm O}(\varepsilon^4),
\monendstar
and we deduce
\moneqstar  \left\{ \begin{array}{l} 
 \Psi_3 = - \partial_{t_1} (\Sigma \, \Psi_2  +  {1\over2} \Psi_2) - {1\over{\tau_0}} \, \partial_{t_2} (\Sigma \, \Psi_1  +  {1\over2} \Psi_1)
 - {1\over{\tau_0^2}} \, \partial_{t_3} \Phi -  {1\over2} \, \partial_{t_1}^2  (\Sigma \, \Psi_1  +  {1\over2} \Psi_1) \\ \qquad \quad 
 -  {1\over{2 \, \tau_0}} \, \partial_{t_1} \, \partial_{t_2}  \Phi(W) 
 -  {1\over{2 \, \tau_0}} \, \partial_{t_2} \, \partial_{t_1}  \Phi(W)  - {1\over6} \, \partial_{t_1}^3 \Phi - D \, (\Sigma \, \Psi_2  -  {1\over2} \Psi_2 )
 \\   \qquad \quad 
 +  {1\over2} \, D_2 \,  (\Sigma \, \Psi_1  -  {1\over2} \Psi_1 ) - {1\over6} \,  ( C_3 \, W  + D_3 \, \Phi ) .
\end{array} \right. \monendstar
We have the auxiliary relations 

\smallskip \noindent
$ - \partial_{t_1} (\Sigma \, \Psi_2  +  {1\over2} \Psi_2) = (\Sigma  +  {1\over2}) \, \dd \Psi_2 . \Gamma_1 $ 

\smallskip \noindent
$  -  {1\over{\tau_0}} \, \partial_{t_2} (\Sigma \, \Psi_1  +  {1\over2} \Psi_1) = (\Sigma  +  {1\over2}) \, \dd \Psi_1 . \Gamma_2 $ 

\smallskip \noindent
$ -  {1\over{\tau_0^2}} \, \partial_{t_3} \Phi = \dd \Phi . \Gamma_3 $

\smallskip \noindent
$ -  {1\over2} \, \partial_{t_1}^2  (\Sigma \, \Psi_1  +  {1\over2} \Psi_1) = 
{1\over2} \,  (\Sigma  +  {1\over2}) \, \partial_{t_1} ( \dd \Psi_1 . \Gamma_1 ) =
- {1\over2} \,  \Sigma  \, \dd ( \dd \Psi_1 . \Gamma_1 ) . \Gamma_1  
-  {1\over4} \, \dd ( \dd \Psi_1 . \Gamma_1 ) . \Gamma_1     $ 

\smallskip  \noindent
$  - {1\over{2 \, \tau_0}} \, \partial_{t_1} \, \partial_{t_2}  \Phi(W) =  {1\over2} \, \partial_{t_1} \, ( \dd \Phi . \Gamma_2 ) $

\smallskip  \noindent
$  - {1\over{2 \, \tau_0}} \, \partial_{t_2} \, \partial_{t_1}  \Phi(W) =  {1\over2} \, \partial_{t_2} \, ( \dd \Phi . \Gamma_1 ) $

\smallskip  \noindent
$ - D \, (\Sigma \, \Psi_2  -  {1\over2} \Psi_2 ) = - D \, \Sigma \, \Psi_2 +  {1\over2} \, D \, \Psi_2 $

\smallskip  \noindent
$  {1\over2} \, D_2 \,  (\Sigma \, \Psi_1  -  {1\over2} \Psi_1 ) =  {1\over2} \, D_2 \, \Sigma \, \Psi_1 -  {1\over4} \, D_2 \, \Psi_1 $

\smallskip  \noindent
and we obtain
\moneq \label{psi3-intermediaire}  \left\{ \begin{array}{l}
 \Psi_3 = (\Sigma  +  {1\over2}) \, \dd \Psi_2 . \Gamma_1 +  (\Sigma  +  {1\over2}) \, \dd \Psi_1 . \Gamma_2 +  \dd \Phi . \Gamma_3
- {1\over2} \,  \Sigma  \, \dd ( \dd \Psi_1 . \Gamma_1 ) . \Gamma_1  \\ \qquad 
-  {1\over4} \, \dd ( \dd \Psi_1 . \Gamma_1 ) . \Gamma_1 
+ {1\over2} \, \partial_{t_1} \, ( \dd \Phi . \Gamma_2 ) + {1\over2} \, \partial_{t_2} \, ( \dd \Phi . \Gamma_1 )  - {1\over6} \, \partial_{t_1}^3 \Phi
- D \, \Sigma \, \Psi_2 +  {1\over2} \, D \, \Psi_2  \\ \qquad \quad 
+ {1\over2} \, D_2 \, \Sigma \, \Psi_1 -  {1\over4} \, D_2 \, \Psi_1 - {1\over6} \,  ( C_3 \, W  + D_3 \, \Phi ) .
\end{array} \right. \monend
We have also the final intermediate relation

\smallskip  \noindent
$ - {1\over6} \, \big[ \partial_{t_1}^3 \Phi + ( C_3 \, W  + D_3 \, \Phi ) \big] = 
 {1\over6} \, \big[  D_2 \, \Psi_1 + D \, \dd \Psi_1 . \Gamma_1 +  \dd ( \dd \Psi_1 . \Gamma_1 ) . \Gamma_1  \big] . $

\smallskip  \noindent
This last relation comes from the following calculation:

\smallskip  \noindent
$ \partial_{t_1} \Phi = - \dd \Phi . \Gamma_1  = - \Psi_1 - C \, W - D \, \Phi $ 

\smallskip
$ \partial_{t_1}^2 \Phi = \dd \Psi . \Gamma_1 + C \, \Gamma_1 + D \, \dd \Phi . \Gamma_1 $

\smallskip \qquad 
$ \, =  \dd \Psi . \Gamma_1 + C \, ( A \, W + B \, \Phi ) + D \, ( \Psi_1 + C \, W + D \, \Phi ) $ 

\smallskip \qquad 
$ \, = C_2 \, W + D_2 \, \Phi + D \, \Psi_1 +  \dd \Psi . \Gamma_1 $
\hfill because 
 $ \,\, C_2 =  C \,  A + D \, C $, $ \, D_2 = C \, B + D^2   $ 

\smallskip
$ \partial_{t_1}^3 \Phi = - C_2 \, \Gamma_1 - D_2 \, ( \Psi_1 + C \, W + D \, \Phi ) - D \, \dd \Psi . \Gamma_1
-  \dd ( \dd \Psi_1 . \Gamma_1 ) . \Gamma_1 $

\smallskip \qquad 
$ \, = - C_2 \, ( A \, W + B \, \Phi)  -  D_2 \, \Psi_1  - D_2 \, (  C \, W + D \, \Phi )
- D \, \dd \Psi . \Gamma_1
-  \dd ( \dd \Psi_1 . \Gamma_1 ) . \Gamma_1 $

\smallskip \qquad 
$ \, = - C_3 \, W - D_3 \, \Phi  -  D_2 \, \Psi_1 - D \, \dd \Psi . \Gamma_1
-  \dd ( \dd \Psi_1 . \Gamma_1 ) . \Gamma_1 $

\hfill because 
$ \,\, C_3 =  C_2 \,  A + D_2 \, C $, $ \, D_3 = C_2 \, B + D_2 \,D   $ 

and $ \,\, -\partial_{t_1}^3 \Phi  - C_3 \, W - D_3 \, \Phi = 
 D_2 \, \Psi_1 + D \, \dd \Psi . \Gamma_1 +  \dd ( \dd \Psi_1 . \Gamma_1 ) . \Gamma_1 $.
 We insert this relation into the expression (\ref{psi3-intermediaire}) and we obtain 

\smallskip \noindent 
$ \Psi_3 = \dd \Phi . \Gamma_3
+ \Sigma \, \dd \Psi_1 . \Gamma_2  +  {1\over2} \, \dd \Psi_1 . \Gamma_2 
+   \Sigma \,\dd \Psi_2 . \Gamma_1  +  {1\over2} \, \dd \Psi_2 . \Gamma_1 - {1\over2} \,  \Sigma \, \dd ( \dd \Psi_1 . \Gamma_1 ) . \Gamma_1 
-  {1\over4} \, \dd ( \dd \Psi_1 . \Gamma_1 ) . \Gamma_1 $

\smallskip \qquad $  
+  {1\over2} \,   \partial_{t_1} \, ( \dd \Phi . \Gamma_2 ) 
+ {1\over2} \,   \partial_{t_2} \, ( \dd \Phi . \Gamma_1 )
- D \, \Sigma \, \Psi_2 +  {1\over2} \, D \, \Psi_2
+ {1\over2} \, D_2 \, \Sigma \, \Psi_1 -  {1\over4} \, D_2 \, \Psi_1 $

\smallskip \qquad \quad $  
+ {1\over6} \, \big[   D_2 \, \Psi_1  + D \, \dd \Psi_1 . \Gamma_1
+ \,  \dd ( \dd \Psi_1 . \Gamma_1 ) . \Gamma_1  \big] $

\smallskip \noindent \quad $ \,\, =
\dd \Phi . \Gamma_3 + \Sigma \, \dd \Psi_1 . \Gamma_2  
+ {1\over2} \, \dd \big(  \dd \Phi . \Gamma_1  - C \, W -  D \, \Phi  ) . \Gamma_2
+  \Sigma \,\dd \Psi_2 . \Gamma_1 
+  {1\over2} \, \dd \big( \Sigma \, \dd \Psi_1  .  \Gamma_1 $

\smallskip \qquad $   + \,   \dd \Phi  .  \Gamma_2 
-    D \, \Sigma \, \Psi_1   \big)
. \Gamma_1 - {1\over2} \,    \Sigma \, \dd ( \dd \Psi_1 . \Gamma_1 ) . \Gamma_1  
 -   {1\over12} \,  \dd ( \dd \Psi_1 . \Gamma_1 ) . \Gamma_1 
+  {1\over2} \,   \partial_{t_1} \, ( \dd \Phi . \Gamma_2 ) 
+ {1\over2} \,   \partial_{t_2} \, ( \dd \Phi . \Gamma_1 ) $

\smallskip \qquad  \quad $  
 - D \, \Sigma \, \Psi_2 +  {1\over2} \, D \,  (  \Sigma \, \dd \Psi_1 .  \Gamma_1 
+    \dd \Phi  .  \Gamma_2  -   D \, \Sigma \, \Psi_1 )
+ {1\over2} \, D_2 \, \Sigma \, \Psi_1
-   {1\over12}  \, D_2 \, \Psi_1
+ {1\over6} \,  D \, \dd \Psi_1 . \Gamma_1  $


\smallskip \noindent \quad $ \,\, =
\dd \Phi . \Gamma_3 + \Sigma \, \dd \Psi_1 . \Gamma_2  
 \,-\, {1\over2} \,  C \, B  \, \Sigma \, \Psi_1 +  \Sigma \,\dd \Psi_2 . \Gamma_1
-  {1\over12} \,  \dd ( \dd \Psi_1 . \Gamma_1 ) . \Gamma_1  
- D \, \Sigma \, \Psi_2   \,-\,  {1\over2} \, D^2  \,  \Sigma \, \Psi_1 $

\smallskip \qquad $
+\, {1\over2} \, D_2  \, \Sigma \, \Psi_1
- {1\over12}  \, D_2  \, \Psi_1 + {1\over6} \,  D \, \dd \Psi_1 . \Gamma_1 $

\smallskip \noindent
after combining several terms. Then we have finally 

\smallskip \noindent
$  \Psi_3  (W) =   \dd \Phi .  \Gamma_3 
+ \Sigma \, \dd \Psi_1  .  \Gamma_2  +    \Sigma \, \dd \Psi_2 .  \Gamma_1
- {1\over12} \, \dd \, (\dd \Psi_1 .  \Gamma_1 ) .  \Gamma_1 
-  D \, \Sigma \, \Psi_2 - {1\over12} \, D_2 \, \Psi_1  
+ {1\over6} \, D \, \dd \Psi_1 .  \Gamma_1  $

\smallskip \noindent
because $ \, D_2 = C \, B + D^2 $.
In this way,  the first relation of (\ref{ordre-4}) is established.
\hfill $\square$

\monitem Chapman-Enskog expansion: Study at order  four

\smallskip \noindent
We look now to the expansion (\ref{developpement-exponentiel}) at order four: 

\smallskip \noindent 
$  m + \varepsilon \, \tau_0 \, \partial_t  m   +   {1\over2} \, \varepsilon^2 \,\tau_0^2 \,    \partial_t^2  m
 +   {1\over6} \, \varepsilon^3 \,\tau_0^3 \,   \partial_t^3  m   +   {1\over24} \, \varepsilon^4 \,\tau_0^4 \,    \partial_t^4  m  $

\smallskip \noindent \qquad $
  = m^*  - \varepsilon \, \tau_0 \,  \Lambda \,  m^*  +   {1\over2} \, \varepsilon^2 \,\tau_0^2 \,  \Lambda^2 \,  m^*
-   {1\over6} \, \varepsilon^3 \,\tau_0^3 \,  \Lambda^3 \,  m^*  + \,  {1\over24} \, \varepsilon^4 \,\tau_0^4 \,  \Lambda^4 \,  m^* +  {\rm O}(\varepsilon^5)  , $

\smallskip \noindent 
and, as usual with the Chapman-Enskog expansion, we replace the time derivative operator
$ \, \partial_t \, $ by 
$ \,  \partial_{t_1} +   \varepsilon \, \partial_{t_2}  +   \varepsilon^2 \, \partial_{t_3}
  +   \varepsilon^3 \, \partial_{t_4} + {\rm O}(\varepsilon^4) $ and deduce
  
\smallskip   \noindent 
$ m + \varepsilon \, \tau_0 \, ( \partial_{t_1} +   \varepsilon \, \partial_{t_2}  +   \varepsilon^2 \, \partial_{t_3}
+   \varepsilon^3 \, \partial_{t_4} + {\rm O}(\varepsilon^4)  )\,  m
+   {1\over2} \, \varepsilon^2 \,\tau_0^2 \,    ( \partial_{t_1}
+  \varepsilon \, \partial_{t_2} +  \varepsilon^2 \, \partial_{t_3} + {\rm O}(\varepsilon^3) )^2 \, m $

\smallskip   \noindent  \quad
$ + \,  {1\over6} \, \varepsilon^3 \,\tau_0^43\,   ( \partial_{t_1} +  \varepsilon \, \partial_{t_2} + {\rm O}(\varepsilon^2) )^3 \, m
+   {1\over24} \, \varepsilon^4 \,\tau_0^4 \,   ( \partial_{t_1} + {\rm O}(\varepsilon) )^4 \, m $

\smallskip   \noindent  \quad \quad $
= m^*  - \varepsilon \,\tau_0 \,   \Lambda \,  m^*
+ {1\over2} \, \varepsilon^2 \,\tau_0^2 \,  \Lambda^2 \,  m^* - {1\over6} \, \varepsilon^3 \,\tau_0^3 \,  \Lambda^3 \,  m^* 
+   {1\over24} \, \varepsilon^4 \,\tau_0^4 \,  \Lambda^4 \,  m^*  +  {\rm O}(\varepsilon^5).  $ 

\smallskip   \noindent
We expand the powers of the noncommutative operators,

\smallskip   \noindent 
$ m + \varepsilon \, \tau_0 \, ( \partial_{t_1} +   \varepsilon \, \partial_{t_2}  +   \varepsilon^2 \, \partial_{t_3}
+   \varepsilon^3 \, \partial_{t_4} + {\rm O}(\varepsilon^4)  )\,  m $

\smallskip   \noindent \quad $
+   {1\over2} \, \varepsilon^2 \,\tau_0^2 \,   [ \partial_{t_1}^2 +  \varepsilon \, ( \partial_{t_1} \, \partial_{t_2} + \partial_{t_2} \, \partial_{t_1} )
+ \varepsilon^2 \, ( \partial_{t_2}^2 +  \partial_{t_1} \, \partial_{t_3}  + \partial_{t_3} \, \partial_{t_1} ) ] \, m  $

 \smallskip  \noindent \qquad  $ +   {1\over6} \, \varepsilon^3 \,\tau_0^3 \, [ \partial_{t_1}^3 \, m  
+  \varepsilon \, (  \partial_{t_2} \, \partial_{t_1}^2 + \partial_{t_1} \, \partial_{t_2} \, \partial_{t_1} + \partial_{t_1}^2 \, \partial_{t_2} ) ] \, m
+   {1\over24} \, \varepsilon^4  \,\tau_0^4 \,   \partial_{t_1}^4  m   $

\smallskip   \noindent  \qquad \quad $
= m^*  - \varepsilon \,\tau_0 \,   \Lambda \,  m^*
+ {1\over2} \, \varepsilon^2 \,\tau_0^2 \,  \Lambda^2 \,  m^* - {1\over6} \, \varepsilon^3 \,\tau_0^3 \,  \Lambda^3 \,  m^* 
+   {1\over24} \, \varepsilon^4 \,\tau_0^4 \,  \Lambda^4 \,  m^*  +  {\rm O}(\varepsilon^5)  $.

\smallskip   \noindent
We take the first component of the above relation,

\smallskip   \noindent 
$ W+ \varepsilon \, \tau_0 \, ( \partial_{t_1} +   \varepsilon \, \partial_{t_2}  +   \varepsilon^2 \, \partial_{t_3}
+   \varepsilon^3 \, \partial_{t_4} + {\rm O}(\varepsilon^4)  )\,  W $

\smallskip   \noindent \quad $
+   {1\over2} \, \varepsilon^2 \,\tau_0^2 \,   [ \partial_{t_1}^2 +  \varepsilon \, ( \partial_{t_1} \, \partial_{t_2} + \partial_{t_2} \, \partial_{t_1} )
+ \varepsilon^2 \, ( \partial_{t_2}^2 +  \partial_{t_1} \, \partial_{t_3}  + \partial_{t_3} \, \partial_{t_1} ) ] \, W  $

 \smallskip  \noindent \qquad  $ +   {1\over6} \, \varepsilon^3 \,\tau_0^3 \, [ \partial_{t_1}^3 \, W  
+  \varepsilon \, (  \partial_{t_2} \, \partial_{t_1}^2 + \partial_{t_1} \, \partial_{t_2} \, \partial_{t_1} + \partial_{t_1}^2 \, \partial_{t_2} ) ] \, m
+   {1\over24} \, \varepsilon^4  \,\tau_0^4 \,   \partial_{t_1}^4  W  = W -  \varepsilon \,\tau_0 \,  ( A \, W + B \,   Y^*  )   $ 

 \smallskip  \noindent \qquad \quad  $
+   {1\over2} \, \varepsilon^2 \,\tau_0^2 \,   ( A_2 \, W + B_2 \,  Y^*  )
- {1\over6} \, \varepsilon^3 \, \tau_0^3 \,   ( A_3 \, W + B_3 \,  Y^*  )   
+ {1\over24} \, \varepsilon^4 \, \tau_0^4 \,   ( A_4 \, W + B_4 \,  Y^*  )  +  {\rm O}(\varepsilon^5)  $

 \smallskip  \noindent
 and we introduce the relation
 \moneqstar
 Y^*  = \Phi(W) + \varepsilon \,\tau_0 \,   \Big( \Sigma \,  -  {1\over2} \Big) \, \Psi_1  
 + \varepsilon^2 \,\tau_0^2 \,   \Big( \Sigma -  {1\over2} \Big)\, \Psi_2   + \varepsilon^3 \,\tau_0^3 \,   \Big( \Sigma -  {1\over2} \Big)\, \Psi_3 
+ {\rm O}(\varepsilon^4)
\monendstar 
into the right-hand side of the previous expression. Then  we obtain the identity

\smallskip   \noindent 
$ W+ \varepsilon \, \tau_0 \, ( \partial_{t_1} +   \varepsilon \, \partial_{t_2}  +   \varepsilon^2 \, \partial_{t_3}
+   \varepsilon^3 \, \partial_{t_4} + {\rm O}(\varepsilon^4)  )\,  W $

\smallskip   \noindent \quad $
+   {1\over2} \, \varepsilon^2 \,\tau_0^2 \,   [ \partial_{t_1}^2 +  \varepsilon \, ( \partial_{t_1} \, \partial_{t_2} + \partial_{t_2} \, \partial_{t_1} )
+ \varepsilon^2 \, ( \partial_{t_2}^2 +  \partial_{t_1} \, \partial_{t_3}  + \partial_{t_3} \, \partial_{t_1} ) ] \, W  $

 \smallskip  \noindent \quad  $ +   {1\over6} \, \varepsilon^3 \,\tau_0^3 \, [ \partial_{t_1}^3 \, W 
+  \varepsilon \, (  \partial_{t_2} \, \partial_{t_1}^2 + \partial_{t_1} \, \partial_{t_2} \, \partial_{t_1} + \partial_{t_1}^2 \, \partial_{t_2} ) ] \, W
+   {1\over24} \, \varepsilon^4  \,\tau_0^4 \,   \partial_{t_1}^4  W  $

 \smallskip  \noindent \quad  $
= W -  \varepsilon \,\tau_0 \, A \, W  
- \varepsilon \,\tau_0 \, B \, [  \Phi(W) + \varepsilon \,\tau_0 \,  (\Sigma \,  -  {1\over2} ) \, \Psi_1  
+ \varepsilon^2 \,\tau_0^2 \,  (\Sigma -  {1\over2} )\, \Psi_2   + \varepsilon^3 \,\tau_0^3 \,  (\Sigma -  {1\over2} )\, \Psi_3 ] $

 \smallskip  \noindent \qquad  $
+   {1\over2} \, \varepsilon^2 \,\tau_0^2 \,  A_2 \, W 
+    {1\over2} \, \varepsilon^2 \,\tau_0^2 \,   B_2 \,[  \Phi(W)
+ \varepsilon \,\tau_0 \,  (\Sigma \,  -  {1\over2} ) \, \Psi_1  + \varepsilon^2 \,\tau_0^2 \,  (\Sigma -  {1\over2} )\, \Psi_2  ] $

 \smallskip  \noindent \qquad  $
- {1\over6} \, \varepsilon^3 \,\tau_0^3 \,   A_3 \, W 
- {1\over6} \, \varepsilon^3 \,\tau_0^3 \, B_3  \,[  \Phi(W) + \varepsilon \,\tau_0 \,  (\Sigma \,  -  {1\over2} ) \, \Psi_1 ]  
+ {1\over24} \, \varepsilon^4 \,\tau_0^4 \,   ( A_4 \, W + B_4 \,  \Phi  )  +  {\rm O}(\varepsilon^5),  $

\smallskip   \noindent
and we identify the fourth-order terms relative to $ \, \varepsilon $:

\smallskip \noindent
$ \partial_{t_4}   W
+    {1\over{2\,\tau_0^2}} \,  (   \partial_{t_2}^2 +  \partial_{t_1} \, \partial_{t_3}  + \partial_{t_3} \, \partial_{t_1}   )  W
+   {1\over{6\,\tau_0}} \, (  \partial_{t_2} \, \partial_{t_1}^2 + \partial_{t_1} \, \partial_{t_2} \, \partial_{t_1}
  + \partial_{t_1}^2 \, \partial_{t_2} )  W  +   {1\over24} \,   \partial_{t_1}^4  W $

\smallskip  \noindent \quad $
+  B  \, (\Sigma -  {1\over2} )\, \Psi_3 
-   {1\over2} \, B_2  \, (\Sigma -  {1\over2} )\, \Psi_2 
+    {1\over6} \, B_3  \,(\Sigma -  {1\over2} )\, \Psi_1
-  {1\over24} \,    ( A_4 \, W + B_4 \, \Phi ) = 0 . $

\smallskip  \noindent
Due to the fourth relation of (\ref{dvpt-dtW-ordre-4}), we can write 

\smallskip \noindent 
$ \Gamma_4 =  B \,  (\Sigma -  {1\over2} )\, \Psi_3 +    {1\over{2\,\tau_0^2}}  \,
 (\partial_{t_2}^2 +  \partial_{t_1}  \partial_{t_3}  + \partial_{t_3}  \partial_{t_1}  )  W
- {1\over2} \, B_2  \,  (\Sigma -  {1\over2} )\, \Psi_2   $           

\smallskip \noindent   \qquad $
+  {1\over{6\,\tau_0}} \, ( \partial_{t_2} \, \partial_{t_1}^2 + \partial_{t_1}  \partial_{t_2}  \partial_{t_1}
+ \partial_{t_1}^2 \, \partial_{t_2} ) W +   {1\over6} \, B_3  \,(\Sigma -  {1\over2} )\, \Psi_1 
+ \,   {1\over24} \, [  \partial_{t_1}^4  W - (  A_4 \, W + B_4 \, \Phi ) ] . $

\smallskip \noindent
We make explicit the following algebraic expressions for intermediate terms:

\smallskip \noindent 
$   {1\over{\tau_0^2}} \, \partial_{t_2}^2  W = -  {1\over{\tau_0}} \, \partial_{t_2} (\Gamma_2)
= - {1\over{\tau_0}} \,\dd \Gamma_2 . \partial_{t_2} W  = \dd \Gamma_2 . \Gamma_2 , $ 

\smallskip \noindent 
$  {1\over{\tau_0^2}} \,  \partial_{t_1}  \partial_{t_3}  W =  - \partial_{t_1} \Gamma_3 =  - \dd \Gamma_3 . \partial_{t_1} W   = \dd \Gamma_3 . \Gamma_1 , $ 

\smallskip \noindent 
$   {1\over{\tau_0^2}} \, \partial_{t_3}  \partial_{t_1} W 
= \dd \Gamma_1 . \Gamma_3 = A \, \Gamma_3 + B \, \dd \Phi . \Gamma_3 , $ 

\smallskip \noindent
$   {1\over{\tau_0}} \, \partial_{t_2} \, \partial_{t_1}^2  W  =    {1\over{\tau_0}} \, \partial_{t_2} (\dd \Gamma_1 . \Gamma_1 )
 =    {1\over{\tau_0}} \,  \partial_{t_2} ( A \, \Gamma_1  +  \dd \Phi . \Gamma_1 )
 =  - A \, \dd \Gamma_1 . \Gamma_2 - B \, \dd \, ( \! \dd \Phi . \Gamma_1 ) . \Gamma_2 ,   $

\smallskip \noindent
$   {1\over{\tau_0}} \, \partial_{t_1} \partial_{t_1}  \partial_{t_2} \,W =  \partial_{t_1} (\dd \Gamma_1 . \Gamma_2) =
 \partial_{t_1} (A \, \Gamma_2 + B \, \dd \Phi . \Gamma_2 )  = - A \, \dd \Gamma_2 . \Gamma_1 - B \, \dd \, ( \! \dd \Phi . \Gamma_2 ) . \Gamma_1 , $  

\smallskip \noindent
$  {1\over{\tau_0}} \, \partial_{t_1}^2 \, \partial_{t_2} W  \!=\! \partial_{t_1} (\dd \Gamma_2 . \Gamma_1) = \partial_{t_1} ( B \, \Sigma \, \dd \Psi_1 . \Gamma_1) 
= -  B \, \Sigma \, \dd \, ( \!\dd \Psi_1 . \Gamma_1)  . \Gamma_1 , $

\smallskip \noindent
and in this way we obtain 
\moneq \label{gamma4-intermediaire}  \left\{ \begin{array}{l}
\Gamma_4 =  B \,  (\Sigma -  {1\over2} )\, \Psi_3 +  {1\over2} \,
[ \dd \Gamma_2 . \Gamma_2 + \dd \Gamma_3 . \Gamma_1 +  A \, \Gamma_3 + B \, \dd \Phi . \Gamma_3 
- B_2  \,  (\Sigma -  {1\over2} )\, \Psi_2 ]
\\ \quad           
-   {1\over6} \, [ A \, \dd \Gamma_1 . \Gamma_2 + B \, \dd \, ( \! \dd \Phi . \Gamma_1 ) . \Gamma_2 
+ A \, \dd \Gamma_2 . \Gamma_1 + B \, \dd \, ( \! \dd \Phi . \Gamma_2 ) . \Gamma_1
\\ \qquad 
+  B \, \Sigma \, \dd \, ( \!\dd \Psi_1 . \Gamma_1)  . \Gamma_1 
-  B_3  \,(\Sigma -  {1\over2} )\, \Psi_1 ]
+ \,   {1\over24} \, [  \partial_{t_1}^4  W - (  A_4 \, W + B_4 \, \Phi ) ] . 
\end{array} \right. \monend
We now establish the identity
$ \,\,   \partial_{t_1}^4  W
- (  A_4 \, W + B_4 \, \Phi ) = B_3 \, \Psi_1 + B_2 \, \dd \Psi_1 . \Gamma_1 + B \,  \dd \, ( \!\dd \Psi_1 . \Gamma_1)  . \Gamma_1 $. 

\smallskip \noindent
We have the following relations:

\smallskip \noindent
$   \partial_{t_1}  W = - \Gamma_1 = - ( A \, W + B \, \Phi ) , $

\smallskip \noindent
$   \partial_{t_1}^2 W = A \, \Gamma_1 + B \, \dd \Phi . \Gamma_1 $ 

\smallskip  \noindent  \qquad  $ \,\,\, =
A \, \Gamma_1 + B \, ( \Psi_1 + C \, W + D \, \Phi) $

\smallskip  \noindent \qquad  $ \,\,\, = 
A \, (A \, W + B \, \Phi ) + B \, \Psi_1 + B \, C \, W + B \, D \, \Phi $

\smallskip  \noindent \qquad  $ \,\,\, = 
A_2 \, W + B_2 \Phi + B \, \Psi_1 $ 

\smallskip \noindent 
$   \partial_{t_1}^3 W = - A_2 \, \Gamma_1 - B_2 \, \dd \Phi . \Gamma_1 - B \, \dd \Psi_1 . \Gamma_1 $

\smallskip  \noindent \qquad  $ \,\,\, =
- A_2 \, ( A \, W + B \, \Phi ) - B_2 \,  ( \Psi_1 + C \, W + D \, \Phi)  - B \, \dd \Psi_1 . \Gamma_1 $

\smallskip  \noindent \qquad  $ \,\,\, =
- A_3 \, W - B_3 \, \Phi  - B_2 \,  \Psi_1 - B \, \dd \Psi_1 . \Gamma_1 $

\smallskip \noindent 
$   \partial_{t_1}^4 W = A_3 \, \Gamma_1 + B_3 \, \dd \Phi . \Gamma_1  + B_2 \,  \dd \Psi_1 . \Gamma_1  + B \,  \dd \, ( \!\dd \Psi_1 . \Gamma_1) . \Gamma_1 $ 

\smallskip  \noindent \qquad  $ \,\,\, =
A_3 \, ( A \, W + B \, \Phi ) + B_3 \,  ( \Psi_1 + C \, W + D \, \Phi)  + B_2 \,  \dd \Psi_1 . \Gamma_1  + B \,  \dd \, ( \!\dd \Psi_1 . \Gamma_1) . \Gamma_1 $

\smallskip  \noindent \qquad  $ \,\,\, =
A_4 \, W +  B_4 \, \Phi +  B_3 \, \Psi_1  + B_2 \,  \dd \Psi_1 . \Gamma_1  + B \,  \dd \, ( \!\dd \Psi_1 . \Gamma_1) . \Gamma_1 $

\smallskip  \noindent  
and
$  \,\,  \partial_{t_1}^4   W
- (  A_4 \, W + B_4 \, \Phi ) = B_3 \, \Psi_1 + B_2 \, \dd \Psi_1 . \Gamma_1 + B \,  \dd \, ( \!\dd \Psi_1 . \Gamma_1)  . \Gamma_1 $. 
We replace this relation in the expression (\ref{gamma4-intermediaire}) to find

\smallskip  \noindent  
$ \Gamma_4 =  B \, \Sigma \, \Psi_3  -  {1\over2} \, B \, \Psi_3 +  {1\over2} \,
  ( \dd \Gamma_2 . \Gamma_2 + \dd \Gamma_3 . \Gamma_1 + A \, \Gamma_3 + B \, \dd \Phi . \Gamma_3 ) 
- {1\over2}  \, B_2  \, \Sigma \Psi_2 + {1\over4} \, B_2 \, \Psi_2 $

\smallskip   \noindent \qquad $
-   {1\over6} \, [ A \, \dd \Gamma_1 . \Gamma_2 +  B \, \dd \, ( \! \dd \Phi . \Gamma_1 ) . \Gamma_2
  +  A \, \dd \Gamma_2 . \Gamma_1  + B \, \dd \, ( \! \dd \Phi . \Gamma_2 ) . \Gamma_1
+  B \, \Sigma \, \dd \, ( \!\dd \Psi_1 . \Gamma_1)  . \Gamma_1 ] $

\smallskip \qquad \quad  $ 
+  {1\over6} \,  B_3  \, \Sigma \, \Psi_1  - {1\over12}  \,  B_3 \, \Psi_1
+  {1\over24} \, [  B_3 \, \Psi_1 + B_2 \, \dd \Psi_1 . \Gamma_1 + B \,  \dd \, ( \!\dd \Psi_1 . \Gamma_1)  . \Gamma_1 ] $

\smallskip   \noindent \quad $  \,\, =
B \, \Sigma \, \Psi_3  -  {1\over2} \, B \, [ \Sigma \, \dd \Psi_1  .  \Gamma_2   +   \dd \Phi .  \Gamma_3 
-  D \, \Sigma \, \Psi_2   + \Sigma \, \dd \Psi_2 .  \Gamma_1
+ {1\over6} \, D \, \dd \Psi_1 .  \Gamma_1  - {1\over12} \,  D_2 \, \Psi_1 $

\smallskip   \noindent \qquad $
- {1\over12} \,  \dd \, (\dd \Psi_1 .  \Gamma_1 ) .  \Gamma_1  ]
+  {1\over2} \,   \dd \Gamma_2 . \Gamma_2
+  {1\over2} \,  \dd \Gamma_3 . \Gamma_1 +  {1\over2} \,  A \, [ 
 B \,  \Sigma \, \Psi_2  + {1\over12}  \,   B_2 \,  \Psi_1  - {1\over6} \,   B \,  \dd \Psi_1 . \Gamma_1 ] 
+  {1\over2} \,  B \, \dd \Phi . \Gamma_3   $

\smallskip   \noindent \qquad \quad $ 
- {1\over2}  \,  B_2  \, \Sigma \Psi_2  +  {1\over4} \, B_2 \, \Psi_2
-   {1\over6} \,  A \, \dd \Gamma_1 . \Gamma_2 -   {1\over6} \,  B \, \dd \, ( \! \dd \Phi . \Gamma_1 ) . \Gamma_2
-   {1\over6} \,  A \, \dd \Gamma_2 . \Gamma_1 
-   {1\over6} \, B \, \dd \, ( \! \dd \Phi . \Gamma_2 ) . \Gamma_1  $

\smallskip   \noindent \qquad \qquad $
-   {1\over6} \,   B \, \Sigma \, \dd \, ( \!\dd \Psi_1 . \Gamma_1)  . \Gamma_1
+ {1\over6} \,  B_3  \, \Sigma \, \Psi_1 - {1\over24} \,    B_3 \, \Psi_1  + {1\over24} \,  B_2 \, \dd \Psi_1 . \Gamma_1 
+  {1\over24} \,   B \,  \dd \, ( \!\dd \Psi_1 . \Gamma_1)  . \Gamma_1    $

\smallskip   \noindent
due to the expressions  (\ref{ordre-3}) for $ \, \Gamma_3 \, $  and  (\ref{ordre-4}) for $ \, \Psi_3 $.
We see that ten terms disappear from the previous expression  because $ \, \Sigma \, \dd \Psi_1 . \Gamma_2 = \dd \Gamma_2 . \Gamma_2 $, and
$ \, BD + AB = B_2 \, $ and $ \, B D_2 + A B_2 = B_3 $. Then, after elementary simplification of some fractions,
we replace $\, \Gamma_3 \, $ in the expression $ \, \dd \Gamma_3 . \Gamma_1 \,$ 
by its expression given by   (\ref{ordre-3}) to obtain

\smallskip  \noindent  
$ \Gamma_4 = 
B \, \Sigma \, \Psi_3  
-  {1\over2} \,  B \,  \Sigma \, \dd \Psi_2 .  \Gamma_1   - {1\over12} \,  B \, D \, \dd \Psi_1 .  \Gamma_1 
+  {1\over2} \, \dd \, [   B \,  \Sigma \, \Psi_2  + {1\over12}  \, B_2 \,  \Psi_1   - {1\over6} \,    B \,  \dd \Psi_1 . \Gamma_1 ] . \Gamma_1 $ 

\smallskip   \noindent \qquad $ 
- {1\over12} \, A \, B \,  \dd \Psi_1 . \Gamma_1   +  {1\over4} \, B_2 \, \Psi_2
-   {1\over6} \,  A \, \dd \Gamma_1 . \Gamma_2  
-   {1\over6} \,  B \, \dd \, ( \! \dd \Phi . \Gamma_1 ) . \Gamma_2
-   {1\over6} \,  A \, \dd \Gamma_2 . \Gamma_1  -   {1\over6} \, B \, \dd \, ( \! \dd \Phi . \Gamma_2 ) . \Gamma_1 $

\smallskip   \noindent \qquad \quad $ 
+  {1\over6} \,  B_3  \, \Sigma \, \Psi_1
-   {1\over6} \,   B \, \Sigma \, \dd \, ( \!\dd \Psi_1 . \Gamma_1)  . \Gamma_1
 + {1\over24} \,    B_2 \, \dd \Psi_1 . \Gamma_1 
+  {1\over12} \,    B \,  \dd \, ( \!\dd \Psi_1 . \Gamma_1)  . \Gamma_1   $.

\smallskip  \noindent  
Four terms clearly vanish  and four others may be eliminated because $ \,  BD + AB = B_2 $. Then 

\smallskip  \noindent  
$ \Gamma_4 = 
B \, \Sigma \, \Psi_3    +  {1\over4} \, B_2 \, \Psi_2
-   {1\over6} \,  A \, \dd \Gamma_1 . \Gamma_2  -   {1\over6} \,  B \, \dd \, ( \! \dd \Phi . \Gamma_1 ) . \Gamma_2
-   {1\over6} \, B \, \dd \, ( \! \dd \Phi . \Gamma_2 ) . \Gamma_1 $ 

\smallskip   \noindent \qquad $ 
+  {1\over6} \,  (A \, B_2 + B \, D_2)   \, \Sigma \, \Psi_1  -   {1\over6} \,   B \, \Sigma \, \dd \, ( \!\dd \Psi_1 . \Gamma_1)  . \Gamma_1  $

\smallskip   \noindent \quad $  \,\, =
B \, \Sigma \, \Psi_3     +  {1\over4} \, B_2 \, \Psi_2 -   {1\over6} \,  A \, ( A \, B \, \Sigma \, \Psi_1  + B \, \dd \Phi . \Gamma_2) 
  -   {1\over6} \,  B \, \dd \, ( \! \dd \Phi . \Gamma_1 ) . \Gamma_2
  -   {1\over6} \,  A \, B \, \Sigma \, \dd \Psi_1 . \Gamma_1 $

\smallskip   \noindent \qquad $ 
-   {1\over6} \, B \, \dd \, ( \! \dd \Phi . \Gamma_2 ) . \Gamma_1 
+  {1\over6} \,    A \, B_2   \, \Sigma \, \Psi_1
+  {1\over6} \,  B \, D_2  \, \Sigma \, \Psi_1 -   {1\over6} \,   B \, \Sigma \, \dd \, ( \!\dd \Psi_1 . \Gamma_1)  . \Gamma_1 $

\smallskip   \noindent \quad $  \,\, =
B \, \Sigma \, \Psi_3     +  {1\over4} \, B_2 \, \Psi_2  +   {1\over6} \,  A \, B \,  D \, \Sigma \, \Psi_1 
 -   {1\over6} \,  A \, B \,  \dd \Phi . \Gamma_2   -   {1\over6} \,  B \, \dd \, ( \! \dd \Phi . \Gamma_1 ) . \Gamma_2 
  -   {1\over6} \,  A \, B \, \Sigma \, \dd \Psi_1 . \Gamma_1 $

\smallskip   \noindent \qquad $
-   {1\over6} \, B \, \dd \, ( \! \dd \Phi . \Gamma_2 ) . \Gamma_1 
+  {1\over6} \,  B \, D_2  \, \Sigma \, \Psi_1
-   {1\over6} \,   B \, \Sigma \, \dd \, ( \!\dd \Psi_1 . \Gamma_1)  . \Gamma_1   $

\smallskip   \noindent \quad $  \,\, =
B \, \Sigma \, \Psi_3     +  {1\over4} \, B_2 \, \Psi_2  -   {1\over6} \, A \, B \, \Psi_2
  -   {1\over6} \,  B \, \dd \, ( \! \dd \Phi . \Gamma_1 ) . \Gamma_2
  -   {1\over6} \, B \, \dd \, ( \! \dd \Phi . \Gamma_2 ) . \Gamma_1 
+  {1\over6} \,  B \, D_2  \, \Sigma \, \Psi_1 $

\smallskip   \noindent \qquad $
-   {1\over6} \,   B \, \Sigma \, \dd \, ( \!\dd \Psi_1 . \Gamma_1)  . \Gamma_1   $

\smallskip \smallskip
because $ \,\,  -\Psi_2   =  D \,  \Sigma \, \Psi_1  -   \dd \Phi . \Gamma_2  - \Sigma \, \dd \Psi_1 . \Gamma_1   $.
Then we have 
\moneqstar  \left\{ \begin{array}{l}
\Gamma_4   (W) =
    B \, \Sigma \, \Psi_3  + {1\over4} \, B_2 \, \Psi_2   +  {1\over6} \, B \, D_2 \, \Sigma \, \Psi_1  -   {1\over6} \, A \, B \, \Psi_2
\\ \qquad \quad 
-  {1\over6} \, B \, \big( \dd \, (\dd \Phi .  \Gamma_1)  . \Gamma_2 
    + \dd \, (\dd \Phi .  \Gamma_2) . \Gamma_1  \big)  - {1\over6} \, B \, \Sigma \,\dd \, (\dd \Psi_1 .  \Gamma_1 ) .  \Gamma_1.
\end{array} \right. \monendstar 
In this way, the second relation of (\ref{ordre-4}) is established.
\hfill $\square $

\bigskip \bigskip    \noindent {\bf \large    8) \quad  Survey and conclusions} 

\smallskip \noindent
In this contribution, we have considered the  exponentiation of differential operators, classical for the  BGK
variant of lattice Boltzmann schemes. We have also used an  exponential  iteration
of the  multi-resolution times lattice Boltzmann schemes.
Then the ``ABCD''  block decomposition of the moment-velocity operator matrix
allows one to formulate in a compact way the asymptotic expansion of the lattice Boltzmann schemes
that give rise to  the equivalent nonlinear partial differential equations  of the conserved moments.
We have calculated the coefficients of the expansion up to    order four, 
with  recursive  formulas containing  less than seven  terms.
To do this, the  intensive use of differential calculus is mandatory, but the calculation
has been systematized and the results have followed in an automatic way.

\smallskip \noindent
We have applied this expansion  at second order for the 
compressible Navier-Stokes equations. We have proposed in~\cite{DL22} 
various lattice Boltzmann schemes in two and three spatial dimensions
with a single particle distribution.
%
The higher-order expressions can be useful to set initial conditions
in simulations, in particular when studying behaviours of a given symmetry.
Third-order precision is also very interesting to avoid some defects
of lattice Boltzmann schemes, as studied in~\cite{LDL23}.

\smallskip \noindent
We hope that this work helps to standardize these kinds of calculations, since most all work in this field
involves expansions to fourth order or less.  At the same time, this work is meant to highlight
the remarkable agreement between two so very different approaches, and to remind practitioners that the problem
of establishing this agreement to all orders is still a very open question.

%
\bigskip \bigskip   \noindent {\bf  \large  Acknowledgments }

\noindent
FD  thanks the Centre National de  la Recherche Scientifique for according a ``Delegation'' 
at the International Research Laboratory 3457   in the ``Centre de Recherches Math\'ematiques''
of the Université de Montr\'eal during the period February-July 2021.
A part of this contribution was done during this period.

\bigskip \bigskip      \noindent {\bf  \large  References }


%


\end{document}